\documentclass[a4paper,10pt]{article}
\usepackage[centertags]{amsmath}
\usepackage{amsfonts}
\usepackage{amssymb}
\usepackage{amsthm}
\usepackage{epsfig}
\usepackage{setspace}
\usepackage{ae}
\usepackage{eucal}
\usepackage[usenames]{color}%
\usepackage{graphicx}

\theoremstyle{plain}
\newtheorem{thm}{Theorem}[section]

\theoremstyle{definition}
\newtheorem{defi}[thm]{Definition}

\theoremstyle{remark}
\newtheorem{exmp}[thm]{Example}
\newtheorem{rem}[thm]{Remark}

\title{Trotter product formulas and global regular upper bounds of the Navier Stokes equation solution}
\author{J\"org Kampen }

\begin{document}

\maketitle

\begin{abstract}
Global upper bounds with respect to regular norms of a controlled incompressible Navier Stokes equation solution with regular data are constructed by an infinite scheme, where we work in bounded ZFC with bounded quantifiers and explicit infinitesimals. Trotter product formula representations of solutions with an infinitesimal error are obtained, which simplify for calculi with explicit infinitesimal and make the spatial effects needed in order to obtain global schemes more transparent. These spatial effects are analyzed in terms of  elliptic integral upper bounds of the nonlinear terms in strong dual Sobolev function spaces, where the viscosity damping term of the Trotter product formula can  set off growth of higher order frequency modes caused by the nonlinear terms at each time step up to first order. The growth of the zero modes is caused by the Burgeres term alone and can be controlled by an external control function or by an auto-control damping term introduced via via time dilatation. In this context spatial derivatives effects of the operator imply that the zero mode contributions of the nonlinear terms are exclusively from the non-zero modes of the modes of the previous infinitesimal time step.  In strong spaces (strong polynomial decay) the damping is strong enough such that an upper bound is preserved for the velocity component functions of the Navier Stokes equation itself. 
\end{abstract}

\section{A global upper bound theorem}
Global regular existence for the incompressible Navier Stokes equation 
\begin{equation}\label{Navleray}
\left\lbrace \begin{array}{ll}
\frac{\partial v_i}{\partial t}-\nu\sum_{j=1}^n \frac{\partial^2 v_i}{\partial x_j^2} 
+\sum_{j=1}^n v_j\frac{\partial v_i}{\partial x_j}=-\nabla_i p,\\
\\
\mbox{div}\mathbf{v}=0,\\
\\
\mathbf{v}(0,.)=\mathbf{h},
\end{array}\right.
\end{equation}
to be solved for $\mathbf{v}=\left(v_1,\cdots ,v_n \right)^T$ on the domain $\left[0,\infty \right)\times {\mathbb T}^n$ follows from the existence of global regular upper bounds. 
Here, $\nabla_ip=p_{,i}$ denotes the spatial derivative of the pressure with respect to the $i$th spatial derivative, and ${\mathbb T}^n$ is a torus of dimension $n$ (equivalent to periodic boundary condition). Pressure can be eliminated by writing the equation in (\ref{Navleray})  in Leray projection form (cf. below). For the initial data vector we sometimes write $\mathbf{h}=\left(h_1,\cdots ,h_n \right)^T$.  
In order to formulate a precise statement concerning the upper bound constants in our theorem we reformulate the Navier Stokes equation in terms of an infinite nonlinear ODE with respect to the modes. Writing the velocity component $v_i=v_i(t,x)$ for fixed $t\geq 0$ in the analytic basis $\left\lbrace \exp\left( \frac{2\pi i\alpha x}{l}\right),~\alpha \in {\mathbb Z}^n\right\rbrace $, i.e., in the form
\begin{equation}
v_i(t,x):=\sum_{\alpha\in {\mathbb Z}^n}v_{i\alpha}(t)\exp{\left( \frac{2\pi i\alpha x}{l}\right) },
\end{equation}
the initial value problem in (\ref{Navleray}) is equivalent an infinite ODE initial value problem for the infinite time dependent vector function of velocity modes $v_{i\alpha},~\alpha\in {\mathbb Z}^n,~1\leq i\leq n$, where
\begin{equation}\label{navode200first}
\begin{array}{ll}
\frac{d v_{i\alpha}}{dt}=\sum_{j=1}^n\nu \left( -\frac{4\pi^2 \alpha_j^2}{l^2}\right)v_{i\alpha}
-\sum_{j=1}^n\sum_{\gamma \in {\mathbb Z}^n}\frac{2\pi i \gamma_j}{l}v_{j(\alpha-\gamma)}v_{i\gamma}\\
\\
+\frac{2\pi i\alpha_i}{l}1_{\left\lbrace \alpha\neq 0\right\rbrace}\frac{\sum_{j,k=1}^n\sum_{\gamma\in {\mathbb Z}^n}4\pi^2 \gamma_j(\alpha_k-\gamma_k)v_{j\gamma}v_{k(\alpha-\gamma)}}{\sum_{i=1}^n4\pi^2\alpha_i^2},
\end{array} 
\end{equation}
for all $1\leq i\leq n$ and where for all $\alpha\in {\mathbb Z}^n$ we have $v_{i\alpha}(0)=h_{i\alpha}$. We denote $\mathbf{v}^F=(v^F_1,\cdots v^F_n)^T$ with $n$ infinite vectors $v^F_i=(v_{i\alpha})_{\alpha \in {\mathbb Z}^n}$. For the regularity of the data we essentially assume that there is a constant $C>0$ and a real number $s>1$ such that for all $1\leq i\leq n$ and all $\alpha\in {\mathbb Z}^n$ we have
\begin{equation}\label{condh}
{\big |}h_{i\alpha}{\big |}\leq \frac{C}{1+|\alpha|^{n+s}}.
\end{equation}
This requirement is closely related to the statement that $h_i\in H^{\frac{n}{2}+s}$ for some $s>1$. In this context, note that for $s>0$ the inequality
\begin{equation}
|h_{i\alpha}|\leq \frac{C}{1+|\alpha|^{n+s}}<\infty
\end{equation}
implies  polynomial decay of order $2n+2s$ of the quadratic modes. Then the equivalence for the dual norm
\begin{equation}\label{hmequiv}
h_i\in h^m\equiv h^m\left({\mathbb Z}^n\right)  \mbox{iff} \sum_{\alpha\in {\mathbb Z}^n}|h_{i\alpha}|^{2m}({1+|\alpha|^{2m}})<\infty
\end{equation}
implies that $h_i\in H^{\frac{n}{2}+s}$.  On the other hand we shall observe that 
at each infinitesimal time step (increment $\delta t$ in time), the upper bound of the nonlinear terms, i.e., the Burgers term and the Leray projection term is given in the form 
\begin{equation}\label{cC0}
2\pi(n+n^2)\sum_{\beta\in {\mathbb Z}^n}\frac{|\beta| C}{1+|\alpha-\beta|^{n+s}}\frac{ C}{1+|\beta|^{n+s}}\delta t\leq \frac{ cC^2}{1+|\alpha|^{n+2s-1}}\delta t,
\end{equation}
where this upper bound can be obtained by comparison with an elliptic integral. Note the mode $\beta$ in the numerator of (\ref{cC0}) which is due to a spatial derivative in the Burgers term and the spatial derivative of the pressure.
Hence it seems that $H^{\frac{n}{2}+1}$ is the critical space for regularity as we have a contractive property of iterated elliptic intergals of for data in $H^{\frac{n}{2}+s}$ for $s>1$. Note that for $s<1$ the relation in (\ref{cC0}) is non-contractive as $n+2s-1<n+s$ for $s<1$, and this means that energy is transported to higher frequency modes upon iterations. Hence if a global regular solution of the Euler equation ($\nu=0$) with data in $H^{\frac{n}{2}+r}$ for some real number $r<1$ is proposed or a solution of the Navier stokes equation with general time dependent force term, then it seems unlikely that such a solution can be correct because data may be constructed such that the solution becomes singular at any short time (due to the exploding relation (\ref{cC0})). Furthermore, it seems that this transport to high frequency cannot be damped by the viscosity term if the initial data are just in $H^1$ (as is sometimes proposed in the case of the Navier stokes equation). For the damping is only effective at a time step if the mode coefficient $v_{i\alpha}$ is not equal to zero, while an increment of the nonlinear terms is a global effect. In order to construct singular solutions we could start with data with zero modes on a set which is not too small, for example the set of modes $M_0:=\lbrace \alpha\in {\mathbb Z}^n||\alpha| \mbox{is even}\rbrace$. Then we can easily construct data in $h_i\in H^1$ with $h_{i\alpha}=0$ for $\alpha\in M_0$ such that the asymptotic upper bound for  modes is $\frac{c_0\delta t}{1+|\alpha|^{1}}$ for some finite constant $c_0>0$, for  $\alpha\in M_0$ and as $|\alpha|$ becomes large after one infinitesimal time step of length $\delta t$ (no stronger decay). There is no viscosity damping in this first time step, and at a second time step the viscosity damping may not be strong enough in order to dampen an increment  of order $cc_0^2(1+|\alpha|^2)$. Indeed $H^1$ seems to critical, as this example shows, i.e., there are data in $H^{1-\epsilon}$ with arbitrary $ \epsilon >0$ such that the viscosity damping cannot be offset the growth of the nonlinear terms after two infinitesimal time steps. We do not work out this here, because the focus of this work is on a global scheme in strong spaces which can be used in order to design algorithms. 

We state the essential result.
\begin{thm}\label{linearboundthm}
Assume  initial  data $h_i\in H^{\frac{n}{2}+s}=H^{\frac{n}{2}+s}\left({\mathbb T}^n\right)$ for some $s>1$ such that for all $1\leq i\leq n$
\begin{equation}
{\big |}h_{i\alpha}{\big |}\leq \frac{C}{1+|\alpha|^{n+s} }
\end{equation}
for some finite constant $C>0$.
Then there exists a global regular upper bound for the solution of the Navier Stokes equation of the form 
\begin{equation}
{\big |}v_i(t,.){\big |}_{H^{n+s}}\leq C^+(1+t),
\end{equation}
for all time $t\in [0,T]$ fr an arbitrary time horizon, where $C^+$ is a finite constant (independent of the time horizon $t>0$ and such that 
\begin{equation}
C^+=RC >1
\end{equation}
with some constants $R,C>0$ such that for all $1\leq i\leq n$ and $\alpha\in {\mathbb Z}^n$ we have
\begin{equation}
R\geq \frac{c^2}{\min\{\nu,1\}}C^2,
\end{equation}
where
\begin{equation}
c=2\pi(n+n^2)\sum_{\beta\in {\mathbb Z}^n}\frac{1}{(1+|\beta|^{\frac{n}{2}+s})^2}.
\end{equation}
\end{thm}

\begin{rem}
The arguments used in the proof of Theorem \ref{linearboundthm} indicate that
\begin{itemize}
 \item[i)] There is a theoretical possibility that for $\nu>0$ the requirements for the data may be weakened to be in $H^{s}$ for lower exponents $s>0$, but the upper bound estimate of the nonlinear growth term at each Euler step, i.e., the upper bound on the right side of the inequality
 \begin{equation}\label{cC1}
2\pi(n+n^2)\sum_{\beta\in {\mathbb Z}^n}\frac{|\beta| C}{1+|\alpha-\beta|^{r}}\frac{ C}{1+|\beta|^{r}}\delta t\leq \frac{ cC^2}{1+|\alpha|^{2r-n-1}}\delta t,
\end{equation}
is not integrable for $r<2.5$ (for any dimension $n$). Hence, in the case of weaker exponents $r<2.5$ the viscosity damping in a global regular scheme would have to be strong enough in order to dampen a strong transport of energy to higher frequencies. The problem here is that the growth of the nonlinear term is global, i.e., the growth of each mode can be influenced by all other modes of the data from the previous time step, while the viscosity damping is local. As an algorithm the scheme proposed looses some feature of stability in any case for data in $H^r$ with $r\leq 2.5$ with $r=2.5$ critical. Therefore we consider $r>2.5$ in this paper in order to have a contractive property of the upper bounds of the nonlinear growth terms as in (\ref{cC0});  
\item[ii)] we sketch an idea below, which indicates that in case $\nu=0$, i.e., in case of the Euler equation, a global scheme of a regular solution branch may be obtained by the method proposed here if the data satisfy $h_i\in H^r$ for some $r>2.5$. However, in this case of $\nu=0$ the arguments below indicate that there are data in Sobolev spaces with exponent $r<2.5$ such that classical solutions, i.e., solutions in $C^1$ cannot exist for these data.
\end{itemize}

\end{rem}

\section{Idea of proof and some alternative schemes}

For $h_i\in H^{\frac{n}{2}+s},~1 \leq i\leq n$ with $s>1$ we assume
\begin{equation}
\forall \alpha\in {\mathbb Z}^n:~{\big |}h_{i\alpha}{\big |}\leq \frac{C}{1+|\alpha|^{n+s}}
\end{equation}
for some constant $C>0$.
We may assume that $h_{i0}=0$ (otherwise consider the related equation for $v_i-h_{i0}$, i.e., shifted by the constant $h_{i0}$). 
From a Trotter product formula in a calculus with explicit infinitesimals we derive an Euler scheme on an infinitesimal time scale. At time step number $m$ the Trotter product formula adds for an infinitesimal time step $\delta t$ a factor
\begin{equation}
\left( \delta_{ij\alpha\beta}\exp\left(-\nu 4\pi^2 \sum_{i=1}^n\alpha_i^2 \delta t \right)\right)  \left( \exp\left( \left( \left( e_{ij\alpha\beta}(m\delta t)\right)_{ij\alpha\beta}\right)\delta t\right)  \right)
\end{equation}
which is 'correct' up to order $\delta t^2$, and where the matrix $\left( e_{ij\alpha\beta}(m\delta t)\right)_{ij\alpha\beta}$ with $1\leq i,j\leq n$ and $\alpha,\beta \in {\mathbb Z}^n$ is related  to the Euler equation terms in (\ref{navode200first}) and defined in (\ref{eulermatrix}) below. This factor applied to the data $\mathbf{v}^F((m-1)\delta t)$ from the previous time step leads to the value $\mathbf{v}^F(m\delta t)$ of the scheme at the next time step. The factor  $\left( \delta_{ij\alpha\beta}\exp\left(-\nu 4\pi^2 \sum_{i=1}^n\alpha_i^2 \delta t \right)\right)$ is the viscosity damping at each time step and the factor $\left( \exp\left( \left( \left( e_{ij\alpha\beta}(m\delta t)\right)_{ij\alpha\beta}\right)\delta t\right)  \right)$ is related to the nonlinear Euler terms, i.e., they correspond to the Burgers term and the Leray projection term. If we replace in the Trotter product formula at each time step the viscosity damping factor by a (smaller or equal) damping factor
\begin{equation}
\left( \delta_{ij\alpha\beta}\exp\left(-\tilde{\nu}(\alpha) 4\pi^2 \sum_{i=1}^n\alpha_i^2 \delta t \right)\right) 
\end{equation}
with
\begin{equation}
0\leq \tilde{\nu}(\alpha)\leq \nu \mbox{ for all }\alpha\in {\mathbb Z}^n,
\end{equation}
then we get a formula which can be still useful in order to construct upper bounds. In this respect note that the first order approximation of the viscosity damping factor
\begin{equation}
1-4\pi^2\nu \sum_{i=1}^n\alpha_i^2 \delta t
\end{equation}
leads to strong damping for higher modes $|\alpha|$, and (depending on a infinitesimal time step size $\delta t$) it can be convenient at least to use a 'dampened' damping in this sense. 
For higher modes $|\alpha|$ the increments of the nonlinear terms at each infinitesimal time step can be estimated by (a first order derivative of) geometric series of elliptic upper bound integrals using the relation
\begin{equation}\label{cC}
2\pi(n+n^2)\sum_{\beta\in {\mathbb Z}^n}\frac{|\beta| C}{1+|\alpha-\beta|^{n+s}}\frac{ C}{1+|\beta|^{n+s}}\leq \frac{ cC^2}{1+|\alpha|^{n+2s-2}}.
\end{equation}
Note that the order of decrease with respect to the modes increases by $n+2s-1-n-s=s-1>0$ for $s>1$.  in this inequality, where $c=c(n)$ is a constant depending only on dimension. The factor $n+n^2$ refers to the number of Euler terms in the equation for each component.
We can use the relation in (\ref{cC}) to get an estimate for an upper bound for the growth caused by the nonlinear terms at each time step. 
The growth for a mode $\alpha$ at time step $l+1$ has the upper bound
\begin{equation}\label{growthv}
\frac{ cC^2}{1+|\alpha|^{n+2s-1}}\delta t-\sum_{j=1}^n\nu \left( 4\pi^2 \alpha_j^2\right)|v_{i\alpha}(l\delta t)|\delta t.
\end{equation}
Note that we may estimate with the modulus of the mode in the last term as the damping term has the 'right' (i.e. opposite) sign of the data mode from the previous time step of the Euler scheme.  Now assume that $s>1$ and $\nu >0$. If a mode $v_{i\alpha}$ in (\ref{growthv}) is not too large, say
\begin{equation}
|v_{i\alpha}(l\delta t)|\leq \frac{\frac{C}{2}}{1+|\alpha|^{n+s}},
\end{equation}
then the upper bound in (\ref{growthv}) that for all $1\leq i\leq n$ and $\alpha\in {\mathbb Z}^n$
\begin{equation}\label{valphal+1}
|v_{i\alpha}((l+1)\delta t)|\leq \frac{C}{1+|\alpha|^{n+s}},
\end{equation}
Otherwise, if $|v_{i\alpha}(l\delta t)|\in\left[ \frac{\frac{C}{2}}{1+|\alpha|^{n+s}},\frac{C}{1+|\alpha|^{n+s}}\right]$ and $|\alpha|>0$ is large enough then we get (\ref{valphal+1}) again (becuase $\nu >0$). Hence the growth is controlled if we can control the lower frequencies. 
The next ideas are then related to the control of these lower modes. First, we still have to take the zero modes into account. However, they satisfy
\begin{equation}\label{navode200first2}
\begin{array}{ll}
v_{i0}(l\delta t)=v_{i0}((l-1)\delta t)
-\sum_{j=1}^n\sum_{\gamma \in {\mathbb Z}^n}2\pi i \gamma_jv_{j(-\gamma)}((l-1)\delta t)v_{i\gamma}((l-1)\delta t)
\end{array} 
\end{equation} 
as the viscosity term and the Leray projection become zero at $|\alpha|=0$. We started with the assumption
\begin{equation}
v_{i0}(0)=h_{i0}=0
\end{equation}
without loss of of generality.
We may then define an external control functions $r_i=( r_{i\alpha})_{\alpha\in {\mathbb Z}^n,~1\leq i\leq n}$ with time dependent modes where $r_{i\alpha}(l\delta t)=0$ for all $\alpha \in {\mathbb Z}^n\setminus \left\lbrace  0\right\rbrace$, all $1\leq i\leq n$ and all step number $l$. For $\alpha =0$ and $1\leq i\leq n$ we define the control function zero mode function $r_{i0}(.)$ inductively via the incremental modes $\delta r_{i0}(.)$. At time zero we define $r_{i0}(0)=0$ for all $\alpha\in {\mathbb Z}^n$ and $1\leq i\leq n$. Having defined $r_{i0}((l-1)\delta t$) for some $l\geq 1$ we define
\begin{equation}\label{control1}
\delta r_{i0}(l\delta t):=-\sum_{j=1}^n\sum_{\gamma \in {\mathbb Z}^n}2\pi i \gamma_jv^r_{j(-\gamma)}((l-1)\delta t)v^r_{i\gamma}((l-1)\delta t),
\end{equation}
and $\delta r_{i\alpha}:=0$ for $\alpha\neq 0$.

A controlled velocity function may be defined inductively, where for all $1\leq i\leq n$ and all $\alpha \in {\mathbb Z}^n$
\begin{equation}\label{control2}
v^r_{i\alpha}(0)=v_{i\alpha}(0)=h_i(0)
\end{equation}
and, for $l\geq 1$ we define for all $1\leq i\leq n$ and all $\alpha \in {\mathbb Z}^n$
\begin{equation}\label{control3}
v_{i\alpha}(l\delta t)=v^r_{i\alpha}(l\delta t)+\delta r (l\delta t).
\end{equation}

We then have $v^{r}_{i0}(l\delta t)=0$ by definition and an upper bound for the controlled function $v^r_{i},~1\leq i\leq n$ with modes $v^r_{i\alpha},~1\leq i\leq n,~\alpha\in {\mathbb Z}^n$ leads immediately to a related upper bound for the original velocity function. Note that a time independent upper bound for $v^{r}_{i\alpha},~1\leq i\leq n,~\alpha\in {\mathbb Z}^n$ leads to an upper bound for the uncontrolled velocity modes which has a linear dependence on time. Hence we may identify $v^r_i$ with $v_i$ or assume that the zero modes are zero without loss of generality. We simply forget about the control function $\delta r$ and may use the symbol $r$ for scaling purposes, i.e. simple linear transformations, in the following.

Scaling is not essential but simplifies the argument (it is more essential in order to design stable algorithms). 
For convenience of the reader we derive the Fourier representation of the scaled velocity function $v^*_i,~1\leq i\leq n$ and pressure $p^*,~1\leq i\leq n$ explicitly, where
\begin{equation}
v_i(t,x)=r^{\lambda}v^*_i(\tau,y)=r^{\lambda}v^*_i(r^{\rho}t,r^{\mu}x),
\end{equation}
and
\begin{equation}
p(t,x)=r^{\lambda}p^*(\tau,y)=r^{\lambda}p^*_i(r^{\rho}t,r^{\mu}x).
\end{equation}
It is intended that $r>1$ and $\lambda,\rho,\mu\in {\mathbb R}$.
Note that
\begin{equation}
v_{i,t}=r^{\lambda}v^*_{i,\tau}(r^{\rho}t,r^{\mu}x)r^{\rho}.
\end{equation}
For the analytic expansion
\begin{equation}
\begin{array}{ll}
v_i(t,x)=\sum_{\alpha \in {\mathbb Z}^n}v_{i\alpha}\exp\left(\frac{2\pi i\alpha x}{l}\right)=r^{\lambda}v^*_i(\tau,y)=r^{\lambda}\sum_{\alpha \in {\mathbb Z}^n}v^*_{i\alpha}\exp\left(\frac{2\pi i\alpha y}{l}\right),\\
\\
p(t,x)=\sum_{\alpha \in {\mathbb Z}^n}p_{\alpha}\exp\left(\frac{2\pi i\alpha x}{l}\right)=r^{\lambda}p^*(\tau,y)=r^{\lambda}\sum_{\alpha \in {\mathbb Z}^n}p^*_{\alpha}\exp\left(\frac{2\pi i\alpha y}{l}\right),
\end{array}
\end{equation}
we get the relation
\begin{equation}
\begin{array}{ll}
v_{i\alpha,t}=r^{\lambda+\rho}v^*_{i,\tau}
\end{array}
\end{equation}
for the time derivatives and
\begin{equation}
\begin{array}{ll}
v_{i,j}(t,x)=\sum_{\alpha \in {\mathbb Z}^n}\frac{2\pi i\alpha_j}{l}v_{i\alpha}\exp\left(\frac{2\pi i\alpha x}{l}\right)\\
\\
=r^{\lambda +\mu}\sum_{\alpha \in {\mathbb Z}^n}\frac{2\pi i\alpha_j}{l}v^*_{i\alpha}\exp\left(\frac{2\pi i\alpha y}{l}\right),\\
\\
p_{,j}(t,x)=\sum_{\alpha \in {\mathbb Z}^n}p_{\alpha}\frac{2\pi i\alpha_j}{l}\exp\left(\frac{2\pi i\alpha x}{l}\right)\\
\\
=r^{\lambda+\mu}\sum_{\alpha \in {\mathbb Z}^n}p^*_{\alpha}\frac{2\pi i\alpha_j}{l}\exp\left(\frac{2\pi i\alpha y}{l}\right),
\end{array}
\end{equation}
for the first order spatial derivatives.
Furthermore, for the second order derivatives of the velocity we have
\begin{equation}
\begin{array}{ll}
v_{i,j,j}(t,x)=\sum_{\alpha \in {\mathbb Z}^n}\frac{-4\pi^2 \alpha^2_j}{l^2}v_{i\alpha}\exp\left(\frac{2\pi i\alpha x}{l}\right)\\
\\
=r^{\lambda +2\mu}\sum_{\alpha \in {\mathbb Z}^n}\frac{-4\pi^2 \alpha^2_j}{l^2}v^*_{i\alpha}\exp\left(\frac{2\pi i\alpha y}{l}\right).
\end{array}
\end{equation}
Hence, the infinite ODE, where for all $1\leq i\leq n$ and for all $\alpha \in {\mathbb Z}^n$
\begin{equation}\label{navode}
\begin{array}{ll}
\frac{d v_{i\alpha}}{dt}=\nu\sum_{j=1}^n \left( -\frac{4\pi \alpha_j^2}{l^2}\right)v_{i\alpha}
-\sum_{j=1}^n\sum_{\gamma \in {\mathbb Z}^n}\frac{2\pi i \gamma_j}{l}v_{j(\alpha-\gamma)}v_{i\gamma}\\
\\
-\frac{2\pi i}{l}\alpha_ip_{\alpha},
\end{array} 
\end{equation}
transforms to
\begin{equation}\label{navode*}
\begin{array}{ll}
r^{\lambda+\rho}\frac{d v^*_{i\alpha}}{d\tau}=r^{\lambda+2\mu}\nu\sum_{j=1}^n \left( -\frac{4\pi \alpha_j^2}{l^2}\right)v^*_{i\alpha}
-r^{2\lambda+\mu}\sum_{j=1}^n\sum_{\gamma \in {\mathbb Z}^n}\frac{2\pi i \gamma_j}{l}v^*_{j(\alpha-\gamma)}v^*_{i\gamma}\\
\\
-\frac{2\pi i}{l}r^{\lambda+\mu}\alpha_ip^*_{\alpha}.
\end{array} 
\end{equation}
Here, note that in the context of the term $2\pi i$ we have $i=\sqrt{-1}$ which should not be confused with the index $1\leq i\leq n$, and  $(\alpha -\gamma)$ denotes the subtraction between multiindices understood componentwise, i.e.,
\begin{equation}
(\alpha -\gamma ) =(\alpha_1-\gamma_1,\cdots,\alpha_n-\gamma_n),
\end{equation}
where brackets are added for notational reasons in order to mark separate multiindices.
Hence,
\begin{equation}\label{navode*}
\begin{array}{ll}
\frac{d v^*_{i\alpha}}{d\tau}=r^{2\mu-\rho}\nu\sum_{j=1}^n \left( -\frac{4\pi \alpha_j^2}{l^2}\right)v^*_{i\alpha}
-r^{\lambda+\mu-\rho}\sum_{j=1}^n\sum_{\gamma \in {\mathbb Z}^n}\frac{2\pi i \gamma_j}{l}v^*_{j(\alpha-\gamma)}v^*_{i\gamma}\\
\\
-\frac{2\pi i}{l}r^{\mu-\rho}\alpha_ip^*_{\alpha},
\end{array} 
\end{equation}
In order to eliminate the scaled pressure modes $p^*_{\alpha}$ using Leray projection of the pressure $p$ we transform the Poisson equation
\begin{equation}\label{poisson1}
\Delta p=-\sum_{j,k}v_{j,k}v_{k,j},
\end{equation}
or the equivalent equation where for all $1\leq i\leq n$ and $\alpha\in {\mathbb Z}^n$
\begin{equation}
p_{\alpha}\sum_{i=1}^n\frac{-4\pi^2\alpha_i^2}{l^2}=\sum_{j,k=1}^n\sum_{\gamma\in {\mathbb Z}^n}\frac{4\pi^2 \gamma_j(\alpha_k-\gamma_k)v_{j\gamma}v_{k(\alpha-\gamma)}}{l^2}.
\end{equation}
We get
\begin{equation}
 r^{\lambda+2\mu}p^*_{\alpha}\sum_{i=1}^n\frac{-4\pi^2\alpha_i^2}{l^2}=r^{2\lambda+2\mu}\sum_{j,k=1}^n\sum_{\gamma\in {\mathbb Z}^n}\frac{4\pi^2 \gamma_j(\alpha_k-\gamma_k)v^*_{j\gamma}v^*_{k(\alpha-\gamma)}}{l^2},
\end{equation}
such that
\begin{equation}\label{palpha}
p_{\alpha}=-1_{\left\lbrace \alpha\neq 0\right\rbrace }r^{\lambda}\frac{\sum_{j,k=1}^n\sum_{\gamma\in {\mathbb Z}^n}4\pi^2 \gamma_j(\alpha_k-\gamma_k)v_{j\gamma}v_{k(\alpha-\gamma)}}{\sum_{i=1}^n4\pi^2\alpha_i^2},
\end{equation}
for $\alpha \neq 0$ becomes
\begin{equation}\label{palpha*}
p^*_{\alpha}=-1_{\left\lbrace \alpha\neq 0\right\rbrace }r^{\lambda}\frac{\sum_{j,k=1}^n\sum_{\gamma\in {\mathbb Z}^n}4\pi^2 \gamma_j(\alpha_k-\gamma_k)v^*_{j\gamma}v^*_{k(\alpha-\gamma)}}{\sum_{i=1}^n4\pi^2\alpha_i^2}
\end{equation}
for $\alpha \neq 0$. Plugging this into (\ref{navode*}) we get
\begin{equation}\label{navode**}
\begin{array}{ll}
\frac{d v^*_{i\alpha}}{d\tau}=r^{2\mu-\rho}\nu\sum_{j=1}^n \left( -\frac{4\pi \alpha_j^2}{l^2}\right)v^*_{i\alpha}
-r^{\lambda+\mu-\rho}\sum_{j=1}^n\sum_{\gamma \in {\mathbb Z}^n}\frac{2\pi i \gamma_j}{l}v^*_{j(\alpha-\gamma)}v^*_{i\gamma}\\
\\
+\frac{2\pi i}{l}r^{\mu-\rho+\lambda}\alpha_i1_{\left\lbrace \alpha\neq 0\right\rbrace }\frac{\sum_{j,k=1}^n\sum_{\gamma\in {\mathbb Z}^n}4\pi^2 \gamma_j(\alpha_k-\gamma_k)v^*_{j\gamma}v^*_{k(\alpha-\gamma)}}{\sum_{i=1}^n4\pi^2\alpha_i^2}
\end{array} 
\end{equation}
The parameter $\lambda$ is for numerical purposes only and may be set to $0$.
Note that negative $\lambda$ lead to an increase of initial data by a factor $r^{-\lambda}>1$. The factor $r^{\lambda}$ of the nonlinear Euler terms compensates then for the additional scaling  of the nonlinear terms. We shall assume $\lambda=0$ henceforth. Let is consider how the magnitude of the Euler terms  behave relative to the magnitude of the viscosity damping term.  
A simple case of scaling which leads to a global upper bound is with respect to original time and just with one large spatial transformation parameter $\mu >0$, i.e., with  
\begin{equation}
\mu >1,~\mbox{ and }\rho=\lambda=0~(\mbox{recall that $r>1$ is assumed)}.
\end{equation}
Note that for such a transformation a global regular upper bound for $v^*_i,~1\leq i\leq n$ implies immediately, that the same global regular upper bound holds also for $v_i,~1\leq i\leq n$.
This choice leads to a viscosity damping with the factor $r^{2\mu-\rho}=r^{2\mu}$ such that for the $m+1$th infinitesimal Euler step we have the viscosity term damping 
\begin{equation}
r^{2\mu}\nu\sum_{j=1}^n \left( -\frac{4\pi \alpha_j^2}{l^2}\right)v^*_{i\alpha}(m\delta t)\delta t,
\end{equation}
and a factor $r^{\mu-\rho+\lambda}=r^{\mu}$ for the nonlinear Euler terms, where the modulus of the Euler term increment has the upper bound  (for size $l\geq 1$)
\begin{equation}
\begin{array}{ll}
{\Big |}-r^{\mu}\sum_{j=1}^n\sum_{\gamma \in {\mathbb Z}^n}\frac{2\pi i \gamma_j}{l}v^*_{j(\alpha-\gamma)}v^*_{i\gamma})m\delta t)\delta t\\
\\
+r^{\mu}2\pi i\alpha_i1_{\left\lbrace \alpha\neq 0\right\rbrace }\frac{\sum_{j,k=1}^n\sum_{\gamma\in {\mathbb Z}^n}4\pi^2 \gamma_j(\alpha_k-\gamma_k)v^*_{j\gamma}(m\delta t)v^*_{k(\alpha-\gamma)}(m\delta t)}{\sum_{i=1}^n4\pi^2\alpha_i^2}\delta t{\Big |}\\
\\
\leq r^{\mu}2\pi (n+n^2)\frac{c(C^*)^2}{1+|\alpha|^{n+s}}\delta t
\end{array} 
\end{equation}
for an inductively assumed upper bound
\begin{equation}\label{upboundv*m}
{\Big |}v^*_{j\gamma}(m\delta t){\Big |}\leq \frac{C^*}{1+|\gamma|^{n+s}}
\end{equation}
for each mode $\gamma\in {\mathbb Z^n}$. If for a mode $\gamma\in {\mathbb Z}^n$ we have ${\Big |}v^*_{j\gamma}(m\delta t){\Big |}\leq \frac{C^*}{2(1+|\gamma|^2)}$ than the upper bound (\ref{upboundv*m}) is preserved for step $m+1$. On the other hand, if for a mode $\gamma\in {\mathbb Z}^n$ we have ${\Big |}v^*_{j\gamma}(m\delta t){\Big |}\geq \frac{C^*}{2(1+|\gamma|^2)}$ then for $\gamma=\alpha\neq 0$ and $l=1$ and
\begin{equation}
r^{\mu}\geq \frac{cC^*2\pi(n+n^2)}{2\nu (1+|\alpha|^{n+s})}
\end{equation}
we have
\begin{equation}
\begin{array}{ll}
r^{\mu}\nu\sum_{j=1}^n {\Big |}\left( -\frac{4\pi \alpha_j^2}{l^2}\right)v^*_{i\alpha}(m\delta t){\Big |}\delta t\geq r^{\mu}\nu\sum_{j=1}^n {\Big |}\left( -\frac{4\pi \alpha_j^2}{l^2}\right)\frac{C^*}{2(1+|\alpha|^{n+s})}{\Big |}\delta t\\
\\
\geq {\Big |}r^{\mu}\nu \delta t{\Big |} \frac{C^*}{2(1+|\alpha|^{n+s})}\delta t \geq 2\pi (n+n^2)\frac{c(C^*)^2}{1+|\alpha|^{n+s}}\delta t,
\end{array}
\end{equation}
and possible growth caused by the Euler terms is offset by the viscosity terms.
Hence the upper bound is preserved for all nonzero modes inductively for the purely spatial scaling given, and  this is sufficient as we may consider controlled equivalent scheme where a purely time-dependent control function forces the zero modes to zero.

Next we consider alternative schemes where we control lower modes introducing an auto-controlled scheme. We consider a time-local transformation where we refer to the transformed function by $u^{lc,t_0}_i,~1\leq i\leq n$ and a time global transformation, where we refer to the transformed function by $u^{gl,t_0}_i,~1\leq i\leq n$. 
First concerning the local-time transformation for time $t_0\geq 0$ and some damping parameter $\theta >0$ we consider the family of comparison functions $u^{lc,t_0}_i,~1\leq i\leq n$ with
\begin{equation}
(1+\theta (t-t_0))u^{lc,t_0}_i(s,.)=v_i(t,.),~s=\frac{t-t_0}{\sqrt{1-(t-t_0)^2}},~t\in [t_0,t_0+1),
\end{equation}
which satisfy the equation
\begin{equation}\label{Navleray2}
\left\lbrace \begin{array}{ll}
\frac{\partial u^{lc,t_0}_i}{\partial s}-\nu \sqrt{1-(t-t_0)^2}^3 \sum_{j=1}^n u^{lc,t_0}_{i,j,j}-\frac{\theta\sqrt{1-(t-t_0)^2}^3}{1+\theta (t-t_0)}u^{lc,t_0}_i\\
\\
+\sqrt{1-(t-t_0)^2}^3(1+\theta (t-t_0))\sum_{j=1}^n u^{lc,t_0}_ju^{lc,t_0}_{i,j}=\\
\\ \sqrt{1-(t-t_0)^2}^3(1+\theta (t-t_0))L\left(\mathbf{ u}^{lc,t_0}\right),\\
\\
\mathbf{u}^{lc,t_0}(0,.)=\mathbf{v}(t_0,.),
\end{array}\right.
\end{equation}
where $L\left(\mathbf{ u}^{lc,t_0}\right)$ denotes the Leray projection term on the torus (which is determined by the Leray projection terms of the modes of course).
In the following we use the abbreviations
\begin{equation}
\mu^{lc,0,t_0}:=\sqrt{1-(t-t_0)^2}^3,~\mu^{lc,1,t_0}(\theta)=\sqrt{1-(t-t_0)^2}^3(1+\theta (t-t_0)),
\end{equation}
if this is convenient in order to avoid a formal length of formulas.
Here the coefficient $\theta$ in the additional potential damping term
\begin{equation}\label{potdamp}
\frac{\theta\sqrt{1-(t-t_0)^2}^3}{1+\theta (t-t_0)}u^{lc,t_0}_i
\end{equation}
can be chosen. Note that we can also choose the time size $\Delta'$ of the subscheme for the function $u^{lc,t_0}_i,~1\leq i\leq n$. In any case we have $(t-t_0)\leq \Delta \in (0,1)$ for a classical subscheme time size $\Delta \in (0,1)$, and, e.g., for $\Delta \theta \leq 1$ and large $\theta$ the potential damping in (\ref{potdamp}) becomes relatively large for large mode values $|u_{i\alpha}(l\delta t)|$. The parameter $\theta$ can be adapted such that at each time step the value $u_{i0}(l\delta t)=0$ is enforced for all time step numbers $l$ (this is an alternative to an external control function for the zero modes). The scheme provides an additional stability by means of the additional potential damping term, but this works only together with strong viscosity damping such that the upper bound for the modes $v_{i\alpha}(l\delta t)$ is preserved at each time steo number $l$. Therefor the scheme is more useful for numerical purposes.  

Second concerning the global-time transformation for time $t_0\geq 0$ and some damping parameter $\theta >0$ we consider the family of comparison functions $u^{gl,t_0}_i,~1\leq i\leq n$ with
\begin{equation}
(1+\theta t)u^{gl,t_0}_i(s,.)=v_i(t,.),~s=\frac{t-t_0}{\sqrt{1-(t-t_0)^2}},~t\in [t_0,t_0+1),
\end{equation}
which satisfy the equation
\begin{equation}\label{Navleray2}
\left\lbrace \begin{array}{ll}
\frac{\partial u^{gl,t_0}_i}{\partial s}-\nu \sqrt{1-(t-t_0)^2}^3 \sum_{j=1}^n u^{gl,t_0}_{i,j,j}-\frac{\theta\sqrt{1-(t-t_0)^2}^3}{1+\theta t)}u^{gl,t_0}_i\\
\\
+\sqrt{1-(t-t_0)^2}^3(1+\theta t)\sum_{j=1}^n u^{gl,t_0}_ju^{gl,t_0}_{i,j}=\\
\\ \sqrt{1-(t-t_0)^2}^3(1+\theta t)L\left(\mathbf{ u}^{gl,t_0}\right),\\
\\
\mathbf{u}^{gl,t_0}(0,.)=\mathbf{v}(t_0,.),
\end{array}\right.
\end{equation}
where we may use similar abbreviations
\begin{equation}
\mu^{gl,0,t_0}:=\sqrt{1-(t-t_0)^2}^3,~\mu^{gl,1,t_0}(\theta)=\sqrt{1-(t-t_0)^2}^3(1+\theta t),
\end{equation}
and where the additional potential damping term
\begin{equation}\label{potdamp}
\frac{\theta\sqrt{1-(t-t_0)^2}^3}{1+\theta t}u^{gl,t_0}_i
\end{equation}
has a smaller effect for larger time $t>0$- this leads to time-dependent estimates. 
This scheme can be used in order to obtain upper bounds for global regular solution branches of the Euler equation.
It is interesting to consider a scaled version of this scheme, i.e., we may consider simple linear variable transformations starting with
\begin{equation}\label{scaled1}
(1+\theta \tau)u^{gl,*,\tau_0}_i(\sigma,.)=v^*_i(\tau,.),~s=\frac{\tau-\tau_0}{\sqrt{1-(\tau-\tau_0)^2}},~\tau\in [\tau_0,\tau_0+1),
\end{equation}
which satisfy the equation
\begin{equation}\label{Navleray2*a}
\left\lbrace \begin{array}{ll}
\frac{\partial u^{gl,*,\tau_0}_i}{\partial s}-\nu r^{2\mu-\rho}\sqrt{1-(\tau-\tau_0)^2}^3 \sum_{j=1}^n u^{gl,*,\tau_0}_{i,j,j}-\frac{\theta\sqrt{1-(\tau-\tau_0)^2}^3}{1+\theta \tau}u^{gl,*,\tau_0}_i\\
\\
+r^{\lambda+\mu-\rho}\sqrt{1-(\tau-\tau_0)^2}^3(1+\theta \tau)\sum_{j=1}^n u^{gl,*,\tau_0}_ju^{gl,*,\tau_0}_{i,j}=\\
\\ r^{\lambda+\mu-\rho}\sqrt{1-(\tau-\tau_0)^2}^3(1+\theta \tau)L\left(\mathbf{ u}^{gl,*,\tau_0}\right),\\
\\
\mathbf{u}^{gl,*,\tau_0}(0,.)=\mathbf{v}^*(\tau_0,.),
\end{array}\right.
\end{equation}
where $L\left(\mathbf{ u}^{gl,*,\tau_0}\right)$ denotes the transformed Leray projection term on the torus.
Note that the potential damping term
\begin{equation}\label{potdamp*b}
\frac{\theta\sqrt{1-(\tau-\tau_0)^2}^3}{1+\theta \tau}u^{gl,*,\tau_0}_i
\end{equation}
is not affected by the scaling, but becomes weaker for large time. For large time a strong the parameter $\theta$ has no essential effect anymore and may be set to $1$. 
In this context it is interesting that we can set upper a global scheme with $\lambda=\rho=0$, $r>1$ and $\mu<0$. However first we write down the scheme for $\lambda=0$ and $\mu,\rho \in {\mathbb R}$ in order to observe  different situations. 
We may consider the latter initial value problem on a time interval $[0,\Delta]$ for $\Delta \in (0,1)$.
For equidistant time discretizations which are related by 
\begin{equation}
\delta \sigma=\frac{1}{\sqrt{1-\Delta^2}}\delta \tau
\end{equation}
and get the recursive relation
\begin{equation}\label{trotterlambdalaa**}
\begin{array}{ll}
u^{gl,*,\tau_0}_{i\alpha}(l\delta \sigma)\doteq  {\Big (}u^{gl,*,\tau_0}_{i\alpha}((l-1)\delta \sigma)\left(1-r^{2\mu-\rho}\mu^{gl,0,\tau_0}\nu 4\pi^2 \sum_{j=1}^n \alpha_j^2\delta \sigma\right)\\
\\
 -r^{\mu-\rho}\mu^{gl,1,\tau_0}(\theta)2\pi i\times\\
 \\
 \times \sum_{j=1}^n\sum_{\gamma \in {\mathbb Z}^n}(\alpha_j-\gamma_j)u^{gl,*,\tau_0}_{i(\alpha-\gamma)}((l-1)\delta \sigma)u^{gl,*,\tau_0}_{j\gamma}((l-1)\delta \sigma)\delta \sigma\\
\\
+ r^{\mu-\rho}\mu^{gl,1,\tau_0}(\theta)2\pi i\alpha_i1_{\left\lbrace \alpha\neq 0\right\rbrace}\times\\
\\
\times \frac{
\sum_{j=1}^n\sum_{\gamma \in {\mathbb Z}^n}\sum_{m=1}^n4\pi^2\gamma_j(\alpha_m-\gamma_m)u^{gl,*,\tau_0}_{m(\alpha-\gamma)}((l-1)\delta \sigma)u^{gl,*,\tau_0}_{j\gamma}((l-1)\delta \sigma)}{\sum_{i=1}^n4\pi^2\alpha_i^2}\delta \sigma{\Big )}\\
\\
-\frac{\theta \sqrt{1-(\tau((l-1)\sigma)-\tau_0)^2}^3}{1+\theta (\tau ((l-1)\sigma)} u^{gl,*,\tau_0}_{i\alpha}((l-1)\delta \sigma)\delta \sigma,
\end{array}
\end{equation}
where we recall that
\begin{equation}
\mu^{gl,0,\tau_0}:=\sqrt{1-(\tau-\tau_0)^2}^3,~\mu^{gl,1,\tau_0}(\theta)=\sqrt{1-(\tau-\tau_0)^2}^3(1+\theta \tau).
\end{equation}
For $\nu >0$ in this local scheme we may choose $\mu=\frac{2}{3}\rho>0$ with $\lambda=0$, $\theta=1$, and $r>1$, a choice, which leads to a coefficient of the Leray projection term and the Burgers term of the form $r^{\mu-\rho+\lambda}=r^{-\frac{1}{3}\rho}$ becomes small for large $\rho$ while the viscosity coefficient is $\nu r^{2\mu}=\nu r^{\frac{4}{3}\rho}$  becomes relatively for large $\rho$. 
For a subscheme of a horizon $\Delta \in (0,1)$ the potential damping term for this choice has the lower bound
\begin{equation}
\frac{ \sqrt{1-(\tau((l-1)\sigma)-\tau_0)^2}^3}{1+ \tau ((l-1)\sigma)}\geq 
\frac{ \sqrt{1-\Delta^2}^3}{1+T},
\end{equation}
where this term is just an additional stabilizer in this case. 

This is different if we choose $r>1$  and $\mu<0$, where $\rho$ may be chosen. In this case the potential damping is relatively small, which is unavoidable in a construction of a global regular scheme which works in the viscosity limit. In this case it is much more difficult to obtain a global scheme. A possibility is the construction via viscosity limits of approximative functions $u^{gl,*,\tau_0,\nu}_i,~1\leq i\leq n$ with smoothed functions in the nonlinear terms, which solve equations of the form
\begin{equation}\label{Navleray2*ab}
\left\lbrace \begin{array}{ll}
\frac{\partial u^{gl,*,\tau_0,\nu}_i}{\partial s}-\nu r^{2\mu-\rho}\sqrt{1-(\tau-\tau_0)^2}^3 \sum_{j=1}^n u^{gl,*,\tau_0,\nu}_{i,j,j}-\frac{\theta\sqrt{1-(\tau-\tau_0)^2}^3}{1+\theta \tau}u^{gl,*,\tau_0}_i\\
\\
+r^{\lambda+\mu-\rho}\sqrt{1-(\tau-\tau_0)^2}^3(1+\theta \tau)\sum_{j=1}^n u^{gl,*,\tau_0,\nu}_ju^{gl,*,\tau_0,\nu}_{i}\ast_{sp}G_{\nu,i}=\\
\\ r^{\lambda+\mu-\rho}\sqrt{1-(\tau-\tau_0)^2}^3(1+\theta \tau)L^{\nu}\left(\mathbf{ u}^{gl,*,\tau_0,\nu}\ast_{sp G_{\nu}}\right),\\
\\
\mathbf{u}^{gl,*,\tau_0}(0,.)=\mathbf{v}^*(\tau_0,.),
\end{array}\right.
\end{equation}
where  $L^{\nu}\left(\mathbf{ u}^{gl,*,\tau_0,\nu}\ast_{sp G_{\nu}}\right))$ is an approximative Leray projection term with first order spatial derivatives $u^{gl,*,\tau_0,\nu}_{i,j}$ replaced by $u^{gl,*,\tau_0,\nu}_{i}\ast_{sp}G_{\nu,i}$.
The resulting upper bound is then a global regular upper bound of a solution bracnh of the Euler equation. This is a different construction, and it is therefor considered elsewhere in detail. 


A note about the difference concerning the construction of regular upper bounds in dual spaces and in classical spaces.
In classical spaces and using a Gaussian fundamental $G_{\nu}$ solution of $G_{\nu,t}-\nu\Delta G_{\nu}=0$ we may use local time solution representations in terms of first order spatial derivatives of the Gaussian of the form
\begin{equation}
\begin{array}{ll}
D^{\alpha}_xv_i(t,.)=D^{\alpha}_xv_i(t_0,.)+\sum_{j=1}^nD^{\beta}_x(v_jv_{i,j})\ast G_{\nu ,l}\\
\\
+D^{\beta}_x\left( K_{n,i}\ast_{sp}\sum_{j,k}v_{j,k}v_{k,j}\right)\ast G_{\nu,l} 
\end{array}
\end{equation}
on a time interval $[t_0,t_0+\Delta]$ for $1\leq |\alpha|\leq m\geq 2$ and where $|\beta|+1=|\alpha|$ with $\beta_l+1=\alpha_l$ and $\beta_k=\alpha_k$ for $l\neq k$. Using the incompressibility condition similar representations can in terms of first order derivatives of the Gaussian can be found for the value function itself. We may then use the antisymmetry of the first order derivatives of the Gaussian with respect to the origin in order to get suitable growth estimates.
However, if we look at spatial transformations $v^*_i(t,y)=v^*_i(t,r^{\mu}x)=v_i(t,x)$ (spatial scaling) then there is an essential difference whether we consider $r^{\mu}$, $r>1$ with $\mu >0$ or with  $\mu<0$. In the former case we have
\begin{equation}
\max_{1\leq i\leq n}\sup_{t\in [t_0,t_0+\Delta]}{\big |}v^*_i(t,.){\big |}_{H^m}\leq C
\Longrightarrow \max_{1\leq i\leq n}\sup_{t\in [t_0,t_0+\Delta]}{\big |}v_i(t,.){\big |}_{H^m}\leq C,
\end{equation}
but if $\mu<0$ then we only have
\begin{equation}
\max_{1\leq i\leq n}\sup_{t\in [t_0,t_0+\Delta]}{\big |}v^*_i(t,.){\big |}_{H^m}\leq C
\Longrightarrow \max_{1\leq i\leq n}\sup_{t\in [t_0,t_0+\Delta]}{\big |}v_i(t,.){\big |}_{H^m}\leq r^{m|\mu|}C.
\end{equation}
It is due to the fact that we can use the effects of convolutions with first order derivatives of the Gaussian that we can construct global regular upper bounds.
The matter is different for dual spaces. Here, we have one Trotter product formula for the value function and this formula represents a regular solution (in $H^{n/2+s}$, $s>1$)if the polynomial decay with respect to the order order of the modes is at least of the form $\frac{1}{1+|\alpha|^{n+s}}$. If this can be proved the existence of regular upper bounds in $H^{n/2+s}$ follows immediately.    

\section{Trotter product formulas}
In the following we describe a proof of the standard theorem of the introduction in a framework of analysis with explicit infinitesimals (such as the functional analytic framework by Connes or the framework of nonstandard analysis). We refer to time intervals $[t_0,t_e]$ for finite times $t_0,t_e\geq 0$ in the usual fashion where we assume that real finite numbers plus multiple infinitesimals are included (in nonstandard terminology we may use hyperfinite discretizations of finite time intervals). The advantage of calculi with explicit infinitesimals is that we do not need a local iteration and that the Trotter product formula has a simpler expression. This way we can observe the damping via the viscosity term in each infinitesimal time step more directly. Similarly for the artificial potential damping terms introduced via a time dilatation transformation. Readers not intimate with calculus with explicit infinitesimals may read the following on an intuitive level (as Leibniz or Diderot would have read it), and then it should be possible for the reader to translate the following into a proof without explicit infinitesimals.

In order to have some terminology (abbreviations) available we use some nonstandard analytic terms in the following. The framework with its definitions can be found in the next section. Note that little of this is actually needed, and we show in section 4 that there is a simple functional analytic proof in ZFC which uses a simple elementary construction of infinitesimals and a classical transfinite induction principle (which is most basic set theory) only.
For arbitrary $t_e>0$, a hyperfinite number $N$ and an infinitesimal $\delta t$ with $N\delta t=t_e$ we define a nonstandard scheme
for time steps $$m\delta t \in \left\lbrace 0,\delta t,2\delta t,\cdots ,(N-1)\delta t, N\delta t=t_e\right\rbrace.$$
We may assume that the letter set is an internal set such that we may apply the internal induction principle and the internal set definition principle. The following argument may be rephrased in other calculi with infinitesimal numbers but it seems that the nonstandard calculus is the most convenient in order to define a transfinite scheme.  
We have the infinitesimal Euler scheme
\begin{equation}\label{navode200first}
\begin{array}{ll}
v_{i\alpha}((m+1)\delta t)=v_{i\alpha}(m\delta t)+\sum_{j=1}^n\nu \left( -\frac{4\pi^2 \alpha_j^2}{l^2}\right)v_{i\alpha}(m\delta t)\delta t\\
\\
-\sum_{j=1}^n\sum_{\gamma \in {\mathbb Z}^n}\frac{2\pi i \gamma_j}{l}v_{j(\alpha-\gamma)}(m\delta t)v_{i\gamma}(m\delta t)\delta t\\
\\
+ \frac{2\pi i\alpha_i}{l}1_{\left\lbrace \alpha\neq 0\right\rbrace}\frac{\sum_{j,k=1}^n\sum_{\gamma\in {\mathbb Z}^n}4\pi^2 \gamma_j(\alpha_k-\gamma_k)v_{j\gamma}(m\delta t)v_{k(\alpha-\gamma)}(m\delta t)}{\sum_{i=1}^n4\pi^2\alpha_i^2}\delta t.
\end{array} 
\end{equation}
The last two terms on the right side of (\ref{navode200first}) correspond to the spatial part of the incompressible Euler equation, where for the sake of simplicity of notation we abbreviate (after same renaming with respect to the Burgers term)
\begin{equation}\label{eulermatrix}
\begin{array}{ll}
e_{ij\alpha\gamma}(m dt)=-\frac{2\pi i (\alpha_j-\gamma_j)}{l}v_{i(\alpha-\gamma)}(m\delta t)\\
\\
+ \frac{2\pi i\alpha_i}{l}1_{\left\lbrace \alpha\neq 0\right\rbrace}
4\pi^2 \frac{\sum_{k=1}^n\gamma_j(\alpha_k-\gamma_k)v_{k(\alpha-\gamma)}(m\delta t)}{\sum_{i=1}^n4\pi^2\alpha_i^2}.
\end{array}
\end{equation}
Note that with this abbreviation (\ref{navode200first}) becomes
\begin{equation}\label{navode200second}
\begin{array}{ll}
v_{i\alpha}((m+1)\delta t)=v_{i\alpha}(m\delta t)+\sum_{j=1}^n\nu \left( -\frac{4\pi^2 \alpha_j^2}{l^2}\right)v_{i\alpha}(m\delta t)\delta t\\
\\
+\sum_{j=1}^n\sum_{\gamma\in {\mathbb Z}^n}e_{ij\alpha\gamma}(m\delta t)v_{j\gamma}(m\delta t)\delta t.
\end{array} 
\end{equation}

Similar schemes can be derived for real mode schemes with $\sin,~\cos$-basis of course, but it is easy to check that for real data $h_i$ the above scheme leads to real solutions, i.e.
\begin{equation}
\forall 1\leq i\leq n~\forall m\geq 1~\forall x:~v_{i}((m\delta t,x)\in {\mathbb R}.
\end{equation}
Here we may assume that the velocity components $v_i$ have their values in the field of standard real numbers where standard parts are taken tacitly if internal counterparts of value functions are considered.     
For the sake of simplicity (and without loss of generality) we consider the case $l=1$ in the following.
In this form the damping effect of the unbounded Laplacian is not obvious. Therefore, we 
derive the Trotter product formula stating that for all $t_e=N_0\delta t$ (where $N_0$ may be a natural number or a hyperfinite number) we have
\begin{equation}\label{aa}
\mathbf{v}^{F}(t_e)\doteq \Pi_{m=0}^{N_0-1}\left( \delta_{ij\alpha\beta}\exp\left(-\nu 4\pi^2 \sum_{i=1}^n\alpha_i^2 \delta t \right)\right)  \left( \exp\left( \left( \left( e_{ij\alpha\beta}\right)_{ij\alpha\beta}(m\delta t)\right)\delta t \right) \right) \mathbf{h}^F, 
\end{equation}
and where $\doteq$ means that the identity holds up to an infinitesimal error. Furthermore, the entries in $(\delta_{ij\alpha\beta})$ are Kronecker-$\delta$s which describe the unit $n{\mathbb Z}^n\times n{\mathbb Z}^n$-matrix. The formula in (\ref{aa}) is easily verified by showing that at each time step $m$ 
\begin{equation}
\left( \delta_{ij\alpha\beta}\exp\left(-\nu 4\pi^2 \sum_{i=1}^n\alpha_i^2 \delta t \right)\right)  \left( \exp\left( \left( \left( e_{ij\alpha\beta}(m\delta t)\right)_{ij\alpha\beta}\right)\delta t\right)  \right)\mathbf{v}^F(m\delta t)
\end{equation}
(as a representation for $\mathbf{v}^F(m\delta t)$) solves the equation (\ref{navode200first}) with an error of order $\delta t^2$.
The use of explicit infinitesimals allow us to have an effective use of first order equality $\doteq$ for arbitrary finite time where on an infinitesimal time level we have simplifications of the formula in (\ref{aa}) in the sense that  
\begin{equation}\label{aabb}
\mathbf{v}^{F}(t_e)\doteq \Pi_{m=0}^{N_0-1}\left( \delta_{ij\alpha\beta}\left(1-\nu 4\pi^2\sum_{i=1}^n\alpha_i^2 \delta t \right)\right)  \left( 1+ \left( \left( e_{ij\alpha\beta}\right)_{ij\alpha\beta}(m\delta t)\right)\delta t \right)  \mathbf{h}^F, 
\end{equation}
is also valid up to order $O(\delta t^2)$ (if $\delta t$ is chosen small enough).

For the logician we note that in a nonstandard framework of enlarged universes we would write down an internal function counterpart of (\ref{aa}) and then show that the formula in (\ref{aa}) holds. 
As we have seen $\nu$ can be made large by transformation. Similarly we can introduce another small parameter factor of the nonlinear terms. Trotter product formulas can be easily adjusted. For example, we may consider the function $w_i,~1\leq i\leq n$ with
\begin{equation}
v_i=\lambda' w_i,~\lambda' w_i(0,.)= h_i,
\end{equation}
(corresponding to $\lambda'=r^{\lambda}$ in our scaling observation above), where we have
\begin{equation}\label{Navleray*}
\left\lbrace \begin{array}{ll}
\frac{\partial w_i}{\partial t}-\nu\sum_{j=1}^n \frac{\partial^2 w_i}{\partial x_j^2} 
+\lambda' \rho\sum_{j=1}^n w_j\frac{\partial w_i}{\partial x_j}=\lambda'L^{\lambda'}\left(\mathbf{ w}^{t_0}\right)\\
\\
\mathbf{w}(0,.)=\lambda'^{-1}\mathbf{h},
\end{array}\right.
\end{equation}
where $L^{\lambda'}\left(\mathbf{ w}^{t_0}\right)$ is a transformed Leray projection term.
This means that we get a small parameter with the nonlinear terms. The Trotter product for $\mathbf{w}=(w_1,\cdots ,w_n)^T$ then becomes
\begin{equation}\label{ww}
\begin{array}{ll}
\mathbf{w}^{F}(t_e)\doteq \Pi_{m=0}^{N_0}\left( \delta_{ij\alpha\beta}\exp\left(-4\pi^2\rho\nu \sum_{i=1}^n\alpha_i^2 \delta t \right)\right)\times\\
\\
\times  \left( \exp\left( \left( \left( e^{\lambda'}_{ij\alpha\beta}\right)_{ij\alpha\beta}(m\delta t)\right)\delta t \right) \right) \lambda^{-1}\mathbf{h}^F, 
\end{array}
\end{equation}
where 
\begin{equation}
\begin{array}{ll}
e^{\lambda'}_{ij\alpha\gamma}(m dt)=-\lambda'\frac{2\pi i (\alpha_j-\gamma_j)}{l}w_{i(\alpha-\gamma)}(m\delta t)\\
\\
+\lambda' \frac{2\pi i\alpha_i}{l}1_{\left\lbrace \alpha\neq 0\right\rbrace}
4\pi^2 \frac{\sum_{k=1}^n\gamma_j(\alpha_k-\gamma_k)w_{k(\alpha-\gamma)}(m\delta t)}{\sum_{i=1}^n4\pi^2\alpha_i^2}.
\end{array}
\end{equation}
It is convenient to give $v_i$ an auto-control, i.e. to compare the value function $v_i,\ 1\leq i\leq n$ for time step interval $[t_0,t_0+1)$ smaller than $1$ with a time dilated function $u^{t_0}_i:[0,1)\times {\mathbb T}^n\rightarrow {\mathbb R},~1\leq i\leq n$, where for $\lambda,\theta >0$ we consider
\begin{equation}
\lambda' (1+\theta t)u^{t_0}_i(s,.)=v_i(t,.),~s=\frac{t-t_0}{\sqrt{1-(t-t_0)^2}},
\end{equation}
for $s\in \left[0,1\right)$, and where it is clearly sufficient to prove the preservation of a upper bound for a time step scheme for the function $u_i,~1\leq i\leq n$. Here the parameter $\lambda' >0$ may be small , where $\lambda'>0$ ensures that the nonlinear terms in the equation for $u^{t_0}_i$ get this additional small parameter, while $\theta>0$ may be large compare to $\lambda$ and is used in order to to control the lower frequency terms such that an upper bound is preserved on the global time level (cf. below). In general we shall have $0< \lambda < \theta$ in order to have a comparatively strong damping term.
The scheme for $u^{t_0}_i,~1\leq i\leq n$ becomes
\begin{equation}\label{navode200secondu}
\begin{array}{ll}
u^{t_0}_{i\alpha}((m+1)\delta t)=u^{t_0}_{i\alpha}(m\delta t)+\mu^{0,t_0}\sum_{j=1}^n\nu \left( -\frac{4\pi^2 \alpha_j^2}{l^2}\right)u_{i\alpha}(m\delta t)\delta t\\
\\
+\sum_{j=1}^n\sum_{\gamma\in {\mathbb Z}^n}e^{u,\lambda'}_{ij\alpha\gamma}(m\delta t)u^{t_0}_{j\gamma}(m\delta t)\delta t.
\end{array} 
\end{equation}
where $\mu^{0,t_0}=\sqrt{1-(.-t_0)^2}^3$ is evaluated at $t_0+m\delta t$ , and
\begin{equation}\label{eujk}
\begin{array}{ll}
e^{u,\lambda' ,t_0}_{ij\alpha\gamma}(m \delta t)=-\lambda'\mu^{1,t_0}\frac{2\pi i (\alpha_j-\gamma_j)}{l}u^{t_0}_{i(\alpha-\gamma)}(m\delta t)\\
\\
+ \lambda'\mu^{1,t_0}\frac{2\pi i\alpha_i}{l}\frac{1_{\left\lbrace \alpha\neq 0\right\rbrace}
\sum_{k=1}^n4\pi^2\gamma_j(\alpha_k-\gamma_k)u^{t_0}_{k(\alpha-\gamma)}(m\delta t)}{\sum_{i=1}^n4\pi^2\alpha_i^2}-\mu^{d,t_0}(m\delta t)\delta_{ij\alpha\gamma}
\end{array}
\end{equation}
along with $\mu^{1,t_0}(t):=(1+\theta t)\sqrt{1-(t-t_0)^2}^3$ and $\mu^{d,t_0}(t):=\frac{\theta\sqrt{1-(t-t_0)^2}^3}{1+\theta t}$.
The last term in (\ref{eujk}) is related to the damping term of the equation for the function $u^{t_0}_i,~1\leq i\leq n$. We summarize the role of the parameters related to damping. For arbitrary time horizon $T>0$ a larger $\theta >0$ causes a larger potential  damping term, where
\begin{equation}\label{dampcoeff}
\frac{\theta\sqrt{1-(t-t_0)^2}^3}{1+\theta (t-t_0)}\geq  \frac{\theta\sqrt{1-\Delta^2}^3}{1+\theta T},
\end{equation}
if we consider the comparison function on the time interval $[0,\Delta]$ for some $0<\Delta <1$ to be chosen later.  Consider the infinitesimal Trotter product formula for this variation. Defining $t^0_e=N_0\delta t+t^0_i$ and for some times $0\leq t^0_i$ and $t^0_e-t^0_i< 1$ along with $t^0_e-t^0_i=N_0\delta t$ for some number $N_0$ (which may be a natural number or a hyperfinite number) we get
\begin{equation}\label{trotterlambda}
\begin{array}{ll}
\mathbf{u}^{F,t_0}(t^0_e)\doteq \Pi_{m=0}^{N_0}\left( \delta_{ij\alpha\beta}\exp\left(-\nu \mu^{t_0}\frac{4\pi^2}{l^2}\sum_{i=1}^n\alpha_i^2 \delta t \right)\right)\times\\
\\
\times \left( \exp\left( \left( \left( e^{u,\lambda',t_0}_{ij\alpha\beta}\right)_{ij\alpha\beta}(m\delta t)\right)\delta t\right)  \right) \mathbf{u}^{F,t_0}(t^0_i).
\end{array}
\end{equation}
\begin{rem}
You may expect that we separate the damping term $-\mu^{d,t_0}(m\delta t)\delta_{ij\alpha\gamma}$ in (\ref{eujk}) in order to have another exponential factor of the damping term in (\ref{trotterlambda}). Indeed this leads to another Trotter product formula, but as we have formulas based on explicit infinitesimals and maintain equality up to infinitesimal error, we may use the simple form in (\ref{trotterlambda}).   
\end{rem}
For the scaled functions define in (\ref{scaled1})
and (\ref{Navleray2*a}) we have
\begin{equation}\label{trotterlambda*}
\begin{array}{ll}
\mathbf{u}^{*,F,t_0}(t^0_e)\doteq \Pi_{m=0}^{N_0}\left( \delta_{ij\alpha\beta}\exp\left(-\nu r^{2\mu-\rho} \mu^{t_0}\frac{4\pi^2}{l^2}\sum_{i=1}^n\alpha_i^2 \delta t \right)\right)\times\\
\\
\times \left( \exp\left( \left( \left( e^{u,\lambda,\mu,\rho,t_0}_{ij\alpha\beta}\right)_{ij\alpha\beta}(m\delta t)\right)\delta t\right)  \right) \mathbf{u}^{*,F,t_0}(t^0_i),
\end{array}
\end{equation}
where
\begin{equation}\label{eujk*}
\begin{array}{ll}
e^{u,\lambda,\mu,\rho ,t_0}_{ij\alpha\gamma}(m \delta t)=-r^{\mu-\rho+\lambda}\mu^{1,t_0}\frac{2\pi i (\alpha_j-\gamma_j)}{l}u^{t_0}_{i(\alpha-\gamma)}(m\delta t)\\
\\
+ r^{\mu-\rho+\lambda}\mu^{1,t_0}\frac{2\pi i\alpha_i}{l}\frac{1_{\left\lbrace \alpha\neq 0\right\rbrace}
\sum_{k=1}^n4\pi^2\gamma_j(\alpha_k-\gamma_k)u^{t_0}_{k(\alpha-\gamma)}(m\delta t)}{\sum_{i=1}^n4\pi^2\alpha_i^2}-\mu^{d,t_0}(m\delta t)\delta_{ij\alpha\gamma}
\end{array}
\end{equation}
Note that the symbols $\mu^{0,t_0},\mu^{1,t_0}$ refer to time-dependent coefficients and are not to  be confused with the scaling parameter $\mu$. Note that in these formulas the spatial transformation expressed by the parmater $r^{\mu}$ has a similar effect as considering the problem on domains of different size (different $l$). 

\section{Background of nonstandard analysis, bounded Zermelo and enlargements of universes}

There has been some criticism of nonstandard analysis, especially by A. Connes, and some criticism of this criticism (cf. \cite{KL}). Our point of view is consistent with the view expressed in \cite{KL}. From a theoretical (logical) point of view nonstandard analysis is equiconsistent with ZFC. Infinitesimals can never be exhibited in space-time, neither the infinitesimals of nonstandard analysis nor the compact operators considered by Connes. Connes' interpretations of infinitesimals as compact operators (cf \cite{C}) leads to noncommuting infinitesimals. However, we do dot need this if we consider schemes for classical equations. So nonstandard analysis is a legitimate possible framework for the scheme considered here, but we emphasize that it can be rephrased in Connes' theory. An elementary introduction into nonstandard analysis may be found in \cite{R}.    
In order to have all terms defined we need to say what an infinitesimal number and a hyperfinite number is in an enlarged universe. We remark that the proof can be rephrased in an argument which works only with basic set theory such as the transfinite induction principle (defined there), and where we define infinitesimal numbers directly in ZFC using non-principal ultrafilters bases on the class of co-finite subsets of the set of natural numbers. No specific features of this construction are needed, and interpretations of infinitesimals as compact operators (as is done by Connes) provide an alternative valid framework for the argument of this paper. First we define the latter term starting with a recapture of bounded Zermelo. Then we define the former terms in the framework of enlarged universes.
Bounded Zermelo is ZFC with quantification restricted to existing sets. The language of bounded Zermelo is a normal set theoretic language, with the exception of the restricted quantification rule, i.e. for any formula $\phi(x)$ with the free variable $x$ and sets $a,b$  
\begin{equation}
\exists x\in a \phi(x), ~ \forall x \in b\phi(x)
\end{equation}
are formulas.
Denote the usual set theoretic language with bounded quantifiers by ${\cal L}_R$. We list  the axioms of restricted ZFC (RZFC in symbols).

\begin{itemize}

\item[(E)] (Extensionality) $y=x$ iff, for all $z$, $z\in x$ iff $z\in y$.

\item[(RC)] (Restricted Comprehension) If $\phi(x)$ is an ${\cal L}_R$-formula (with quantifiers restricted)
and free variable $x$ and $a$ is a set, then there exists a set $b$ with $x\in b$ iff $x\in a$ and $\phi(x)$.

\item[(NS)] (Null Set) There exists a set $\oslash$ such that for all $x$ $x\notin \oslash$.

\item[(P)] (Pair) For all $x$ and $y$ there exists $z$ with $u\in z$ iff $u=x$ or $u=y$.

\item[(U)] (Union) For all $x$ there exists $y$ with $z\in y$ iff there exists $w$ with $z\in w\in y$.

\item[(PS)] (Power Set) for all $x$ there exists $y$ with $z\in y$ iff $z\subseteq x$. The set $y$ will be denoted
by $P(x)$ 

\item[(F)] (Foundation) For all $x\neq \oslash $ there exists a set $y\in x$ with $y\cap x=\oslash$.

\item[(I)] (Axiom of Infinity) There exists a set ${\mathbb N}$ such that $\oslash \in {\mathbb N}$ and $x\in {\mathbb N}$ implies  $x\cup \{x\}\in {\mathbb N}$ 

\item[(AC)] (Axiom of Choice) if $I\neq \oslash$ is an index set and for all $i\in I$ $X_i\neq \oslash$, then
$\Pi_{i\in I}X_i\neq \oslash$.

\end{itemize}
Bounded Zermelo is well known to be equiconsistent to a model of the first order theory of well-pointed topoi. However, as explained before, we are not dealing with the topos-theoretic view here. 


Next we consider enlarged universes. Several types of extensions of universes may be considered (systems of sets with certain closure properties). A standard extension is called enlargement of universes is an axiomatic enlargement of an universe $U$ and is defined as an embedding
\begin{equation}
U\stackrel{*}{\rightarrow} U'
\end{equation}
which satisfies some axioms. Another one is by an ultrafilter construction with an index set $I$. So we seek for some object of the form
\begin{equation}
"U^*:=U^I/{\cal F}",
\end{equation}
where ${\cal F}$ is a nonprincipal ultrafilter. Let us sketch this.

The axiomatic system RZFC tells us what sets are. In set theory sometimes you take entities for granted (as Kronecker expresses that the natural numbers are given by god, and the rest is constructed by mankind). Similar in nonstandard theory we tend to consider some set as the set of urelements. In ZFC people realised that urelements are superfluous. However, when we talk about universes over ..., we include urelements as convenient. We denote sets by capital letters and elements which are either sets or urelements by small letters. E.g. in $A\in {\mathbb U}$ $A$ is a set, while in $a\in {\mathbb U}$ a may be a set or an urelement. A universe is a set with certain properties. First if $A\in {\mathbb U}$ is a set, then we want all elements of $a$ to be present in ${\mathbb U}$, i.e.
\begin{equation}\label{U1}
a\in A\in {\mathbb U}~~\Rightarrow~~a\in {\mathbb U}. 
\end{equation}
Any set ${\mathbb U}$ which satisfies \ref{U1} is called transitive. Furthermore, in a universe we want to have with a set $A\in {\mathbb U}$ its transitive closure $\mbox{Tr}(A)\in {\mathbb U}$. Here $\mbox{Tr}(A)$ is the smallest transitive set that contains $A$. If $A$ is transitive itself, then $A=\mbox{Tr}(A)$, of course. We require:
\begin{equation}\label{U2}
\mbox{if $A\in {\mathbb U}$, then there ex. a transitive set $B\in {\mathbb U}$
 with $A\subset B\subset {\mathbb U}$}.
\end{equation}
Finally we require (and are allowed to by RZFC) that
\begin{equation}
\begin{array}{ll}
\mbox{ if $a,b\in U$, then $\{a,b\}\in {\mathbb U}$}\\
\\
\mbox{ if $A,B\in {\mathbb U}$ are sets, then $A\cup B\in {\mathbb U}$}\\
\\
\mbox{ if $A\in {\mathbb U}$ is a set, then $P(A)\in {\mathbb U}$}
\end{array}
\end{equation}
If the universe contains a set ${\mathbb S}$ such that the members of $S$ are individuals in the sense that
\begin{equation}
\forall x \in {\mathbb S}\left[x\neq \oslash \wedge (\forall y\in {\mathbb U}(y\notin x)\right]. 
\end{equation}
The next step of course is to show that universes exist. They are realized by superstructures.
Let ${\mathbb S}$ be a set. We define a series cumulative power set by
\begin{equation}
\begin{array}{ll}
{\mathbb U}_0\left( {\mathbb S}\right)= {\mathbb S},\\
\\
{\mathbb U}_{n+1}\left( {\mathbb S}\right)= {\mathbb U}_{n}\left( {\mathbb S}\right)\cup P\left( {\mathbb U}_{n}\left( {\mathbb S}\right)\right). 
\end{array}
\end{equation}
Then it is easy to check that
\begin{equation}
{\mathbb U}\left( {\mathbb S}\right):=\bigcup_{n\in {\mathbb N}}{\mathbb U}_{n}\left( {\mathbb S}\right)
\end{equation}
is a universe. If the set ${\mathbb S}$ is known from the context, then we . The language ${\cal L}_R$ with quantification restricted to the sets of the universe ${\mathbb U}$ is denoted ${\cal L}_R^{\mathbb U}$. Next, a nonstandard framework for a set ${\mathbb S}$ comprises a universe ${\mathbb U}$ over ${\mathbb S}$ and a map
\begin{equation}
{\mathbb U} \stackrel{*}{\rightarrow}{\mathbb U}', 
\end{equation}
which satisfies
\begin{equation}
\begin{array}{ll}
a^*=a \mbox{ for $a\in {\mathbb S}$},\\
\\
\oslash^*=\oslash,\\
\\
\mbox{ the ${\cal L}_R^{\mathbb U}$-sentence $\phi$ is true iff $\phi^*$ is true}.
\end{array}
\end{equation}  
Such a nonstandard framework is called an enlargement if the following condition is satisfied:
\begin{equation}
\begin{array}{llc}
\mbox{if $A\in {\mathbb U}$ is a collection of sets with the finite intersection property, then there}\\
\mbox{exists an element $z\in {\mathbb U}'$ such that}\\
\\
\hspace{4.5cm}z\in \bigcap\{Z^*| Z\in A \}. 
\end{array}
\end{equation}
\begin{exmp}If ${\mathbb U}$ is a universe on ${\mathbb R}$ and $A$ is the set of intervals $\{(0,r)|r>0\}$, then the latter set satisfies the finite intersection property. Then the enlargement principle tells us that there exists a positive infinitesimals, i.e. there exists
\begin{equation}
b\in \cap \{ (0,r)^*|r>0\}=\mbox{'set of positive infinitesimals'}.
\end{equation}
\end{exmp}

\begin{exmp} Consider a universe $U$ on ${\mathbb N}$ and an enlargement 
$U\stackrel{*}{\rightarrow} U'$.
The set 
\begin{equation}
A=\{{\mathbb N}_{\geq n}| n\in {\mathbb N}\}
\end{equation}
with ${\mathbb N}_{\geq n}:=\{m\in {\mathbb N}| m\geq n\}$
satisfies the finite intersection property. Then the enlargement principle tells us that
\begin{equation}
\exists b: b\in \cap\{ {\mathbb N}_{\geq n}^*|n \in {\mathbb N}\}={\mathbb N}^*\setminus{\mathbb N},
\end{equation}
i.e. there is a set of unlimited numbers.
\end{exmp}

The question now is whether enlargements really exist. This can be shown with the ultrafilter construction. Two types of sets are of special interest for us:
the first type are the internal sets:
\begin{equation}
\mbox{ $a\in {\mathbb U}$ is internal if $a\in A^*$ for some $A\in {\mathbb U}$};
\end{equation}
the second type are the hyperfinite sets: let
\begin{equation}
P_F(A)=\{B\subseteq A| B \mbox{ is finite }\}.
\end{equation}
Then the hyperfinite sets are the members of $P_F(A)^*\in {\mathbb U}'$. 
    

Let $I$ be an infinite set and let ${\cal F}$ be a nonprincipal ultrafilter on $I$. Let ${\mathbb S}$ be a set and let let ${\mathbb U}$ be a universe over ${\mathbb S}$. To $a\in {\mathbb U}$ assign $a_I \in {\mathbb U}^I$, the function with constant value $a$. This way we embed ${\mathbb U}$ in a larger universe.
Similarly as in the ultrafilter construction on ${\mathbb R}$ we have to consider equivalence of elements with respect to the ultrafilter, this time of functions $f,g\in {\mathbb U}^I$. We say
\begin{equation}
\begin{array}{ll}
f\sim g \mbox{ iff } \left\lbrace i|f(i)=g(i)\right\rbrace \in {\cal F}\\
\\
f\in g \sim f' \in g' \mbox{ iff } \left\lbrace i|f(i) \in g(i) \& f'(i) \in g'(i)\right\rbrace \in {\cal F}
\end{array}
\end{equation}
We denote equivalence classes b $\left[ .\right] $ as before and define 
\begin{equation}
\begin{array}{ll}
{\cal W}_n:=\left\lbrace f\in {\mathbb U}_n^I| \left\lbrace i| f(i)\in {\mathbb U}_n\right\rbrace \in {\cal F}\right\rbrace .
\end{array}
\end{equation}
\begin{equation}
\begin{array}{ll}
{\cal W}:=\cup_{n\in {\mathbb N}} {\cal W}_n
\end{array}
\end{equation}
Now for $f\in {\cal W}_0$ let $[f]:=\left\lbrace [h]|h\in {\cal W}_0| \right\rbrace$, and let
\begin{equation}
{\mathbb Y}=\left\lbrace [f]|f\in {\cal W}_0\right\rbrace 
\end{equation}
This defines 
\begin{equation}
{\mathbb U}_0\left({\mathbb Y} \right)={\mathbb Y} 
\end{equation}
Inductively, having defined ${\mathbb U}_n\left({\mathbb Y} \right) $, for $f\in {\cal W}_{n+1}\setminus {\cal W}_n$ define
\begin{equation}
[f]=\left\lbrace [h]|h\in {\cal W}_n \mbox{ and } \left\lbrace i|h(i)\in f(i) \right\rbrace \in {\cal F} \right\rbrace  
\end{equation}
Then
\begin{equation}
{\mathbb U}_{n+1}\left({\mathbb Y} \right)={\mathbb U}_{n}\left({\mathbb Y} \right)\cup\left\lbrace [f]| f\in {\cal W}_{n+1}\setminus {\cal W}_n\right\rbrace ,
\end{equation}
and
\begin{equation}
{\mathbb U}\left({\mathbb Y} \right)=\cup_{n\in {\mathbb N}}{\mathbb U}_{n}\left({\mathbb Y} \right).
\end{equation}
Now, ${\mathbb U}\left({\mathbb Y} \right)$ is the ultrafilter-enlargement we had looked for.
For each $f,g\in {\cal W}$ it is easy to see that
\begin{equation}
\begin{array}{ll}
 [f]\in [g] \mbox{ iff } \left\lbrace i|f(i)\in g(i)\right\rbrace \in {\cal F},~
 [f]= [g] \mbox{ iff } \left\lbrace i|f(i)= g(i)\right\rbrace \in {\cal F}
\end{array} 
 \end{equation}
The map
\begin{equation}
\begin{array}{ll}
*: {\mathbb U}({\mathbb X})\rightarrow  {\mathbb U}({\mathbb Y}),~
  a \rightarrow  a^*=\left[a_I \right]
\end{array} 
\end{equation}
is an embedding of the universe ${\mathbb U}({\mathbb X})$ in the universe ${\mathbb U}({\mathbb Y})$, and we observe that
\begin{equation}
\oslash^*=\oslash, \mbox{ and } {\mathbb X}^*={\mathbb Y}. 
\end{equation}
Since ${\cal F}$ is an ultrafilter, we know that the enlargement ${\mathbb U}({\mathbb Y})$ has nonstandard members.
Let $L_{{\mathbb U}({\mathbb X})}$ and $L_{{\mathbb U}({\mathbb Y})}$ be the formal languages of the respective universes. Denote the model of the ultrafilter-enlargement by ${\cal U}_{{\mathbb Y}}=\left( {\mathbb U}\left({\mathbb Y} \right), \in \right) $ the model of the original universe by ${\cal U}_{{\mathbb X}}=\left( {\mathbb U}\left({\mathbb X} \right), \in \right) $. Then we get the following version of the theorem of Loos.
\begin{thm}
For any $L_{{\mathbb U}({\mathbb X})}$-formula $\phi(x_1,\cdots,x_m)$ and $f_1,\cdots f_m\in {\cal W}$
\begin{equation}
{\cal U}_{{\mathbb Y}}\models \phi \left( e{\big |}^{x_1,\cdots ,x_m}_{[f_1],\cdots ,[f_m]}\right)
\mbox{iff} \left\lbrace i| {\cal U}_{{\mathbb X }}\models\phi\left(e^i{\big |}^{x_1,\cdots ,x_m}_{ f_1(i),\cdots , f_m(i)} \right)  \right\rbrace \in {\cal F}
\end{equation}
\end{thm}
Let $U$ be a universe (which contains the real numbers as individuals) and let $U\stackrel{*}{\rightarrow}U'$ be an enlargement. For $A\in U$ and let
\begin{equation}
P_F(A)=\left\lbrace b\subseteq A| B \mbox{ is finite }\right\rbrace. 
\end{equation}
$P_F(A)^*$ are called hyperfinite subsets of $A$. As an example, consider $P_F({\mathbb N}$ and the $L_U$-sentence
\begin{equation}
\forall n\in {\mathbb N} \exists A \in P_F({\mathbb N}) \forall m \in {\mathbb N} \left[m\in A \leftrightarrow m\leq n \right], 
\end{equation}
i.e. the sentence which has the meaning that for each natural number $n\in {\mathbb N}$ there is a set
$A=\{1,\cdots, n\}$ in 
$P_F({\mathbb N})$. The transfer sentence is
\begin{equation}
\forall n\in {\mathbb N}^* \exists A  \in P_F({\mathbb N})^* \forall m \in {\mathbb N}^* \left[m\in A \leftrightarrow m\leq n \right]. 
\end{equation}
Hence, for all $n\in {\mathbb N}^*$
\begin{equation}
A=\{1,\cdots , n\}\in P_F({\mathbb N})^*. 
\end{equation}
Note that $n$ can be infinite, i.e. $n\in {\mathbb N}^*\setminus {\mathbb N}$.
We prove
\begin{thm}\label{hypfinex}
$A$ is hyperfinite iff there exists $n\in {\mathbb N}^*$ and an internal bijection
\begin{equation}
f:\{1,\cdots ,n\}\rightarrow A.
\end{equation}
\end{thm}
Here a function $f:A\rightarrow B$ is called internal if the set $\mbox{graph}(f)\subseteq A\times B$
is internal.

Proof. Consider a $L_{{\mathbb U}}$-formula
\begin{equation}
\phi\left(X,Y,n,f \right)
\end{equation}
which expresses that $f: X\rightarrow Y$ with
$X=\left\lbrace m\in {\mathbb N}|m\leq n \right\rbrace $
is a bijection. 
Then the $L_{{\mathbb U}}$-sentence
 \begin{equation}
 \psi\equiv \forall Y\in P_F(B) \exists n\in {\mathbb N} \exists f\in P({\mathbb N}\times B)\exists X\in
 P({\mathbb N})\phi\left( X,Y,n,f\right) 
\end{equation}
asserts that for all $Y\in P_F(B)$ there is a number $n$ and a bijection between $X=\left\lbrace 1,\cdots ,n\right\rbrace$  and $Y$- a sentence which is true. The sentence $\psi^*$ is true by transfer. So if $B\in {\mathbb U}$ and  
 $A\in P_F(B)^*$ then the claim follows from the truth of $\psi^*$.
For the converse suppose that there is an internal bijection $f:X=\left\lbrace 1,\cdots n\right\rbrace \rightarrow A$ for some $n\in {\mathbb N}^*$. Then $A$ is internal, because it is the range of an internal function. We want to show that $A$ is hyperfinite. First we observe that
\begin{equation}
\exists X\in P\left({\mathbb N} \right)^* \phi(X,A,n,f)
\end{equation}
is true. Hence the claim that $A$ is hyperfinite follows from transfer of the true $L_{{\mathbb U}}$-sentence
\begin{equation}
\forall Y\in B \exists n\in {\mathbb N} \exists f\in P\left( {\mathbb N}\times A\right) \exists X\in P\left({\mathbb N} \right)\left(  \phi(X,A,n,f)\rightarrow Y\in P_F(A)\right) .
\end{equation}
Having defined enlargements and hyperfinite numbers we now now can easily define all other terms which are used in a nonstandard form of the argument. These are the internal induction principle and the internal set definition principle. For a hyperfinite number the discretization
\begin{equation}
\left\lbrace \frac{k}{N}| k\mbox{ hypernatural }\& k\leq N \right\rbrace=\left\lbrace  0,\frac{1}{N},\frac{2}{N},\cdots ,\frac{N-1}{N},1\right\rbrace 
\end{equation}
is an internal set. Internal sets have the advantage that they obey an internal induction principle based on the fact that each internal subset of the set of hypernaturals has a least member. We have
\begin{thm}\label{intind}
An internal subset $S$ of the set of hypernatural numbers ${\mathbb N}^*$ which contains $1$, i.e., $1\in S$, and is closed under the successor operation $n\rightarrow n+1$, i.e., $n\in S\rightarrow n+1\in S$ equals the whole set of hypernatural numbers, i.e., $S={\mathbb N}^*$.
\end{thm}
Next we recall the internal set definition principle.
\begin{thm}
Let ${\mathbb U}'$ be an enlargement of an universe and let $\psi(x)$ be an internal $L_{ U'}$ formula, where $x$ is the only free variable. Then for any internal set $S\in {\mathbb U}'$ the subset
\begin{equation}
R:=\left\lbrace x\in S|\phi(x)\right\rbrace 
\end{equation}
is internal. Here '$\psi$ is an internal formulas' means that $\psi(x)\equiv \psi(x,a_1,\cdots,a_m)$ for some   internal constants $a_i\in {\mathbb U}'$.
\end{thm}

\section{Proof of the theorem \ref{linearboundthm}}
Let $\nu>0$ and consider data $h_i,~1\leq i\leq n$ with $h_i\in H^{\frac{n}{2}+s}$ for $s>1$ such that for some $C>0$ we have for all $1\leq i\leq n$ and all $\alpha \in {\mathbb Z}^n$
\begin{equation}
|h_{i\alpha}|\leq \frac{C}{1+|\alpha|^{n+s}}.
\end{equation}

 We consider a scheme with viscosity damping for simple transformed velocity component functions $v^*_i,~1\leq i\leq n$, where $r>1$ and $\mu >0$. In this situation we can set up a scheme with $\rho=\lambda=0$. Recall that
\begin{equation}
\max_{1\leq i\leq n}\sup_{t\in [0,T]}{\big |}v^*_i(\tau,.){\big |}_{H^{\frac{n}{2}+s}}\leq C\Longrightarrow \max_{1\leq i\leq n}\sup_{t\in [0,T]}{\big |}v_i(t,.){\big |}_{H^{\frac{n}{2}+s}}\leq C.
\end{equation}
Here we assume that the zero modes are forced to be zero by an external control function.
Note that the linear dependence on the time horizon of the upper bound in the statement of theorem \ref{linearboundthm} is due to an upper bound for the external control function fr the zero modes which has to be added eventually.
Here, recall the reduction to a problem with zero zero modes in section 2, i.e., we assume that the modes $v_{i\alpha},~1\leq i\leq n,~\alpha\in {\mathbb Z}^n$ of velocity component functions $v_i$ satisfy $v_{i0}=0$ for all $1\leq i\leq n$ and $\alpha=0$. As we explained in section 2 this can always achieved by adding an external control function and consider $v^r_i=v_i+r_i$ with the control function $r=(r_1,\cdots,r_n)$ defined as in (\ref{control1}), (\ref{control2}), (\ref{control3}) above.

Next we provide the details of the argument outlined in section 2.
We first consider induction principles which can be used. The following argument works for only for $\nu >0$ because we need strong viscosity damping. In case of $\nu=0$ auto-controlled schemes may be considered.
 Note that auto-controlled schemes can be applied in order to prove the existence of global solution branches but not uniqueness. Especially, they can also be used in order to prove singular solutions. In a situation where there is a singular solution there is often a global solution branch as well. In general it is difficult to prove uniqueness in a situation where we have singular solutions. For example in a related paper we prove the existence of singular solutions for Navier Stokes equation with $L^2$-force data. However, we cannot conclude from this that there is no global strong solution for such models unless we have a proof of uniqueness.
The theoretical minimum for the proof is a calculus with infinitesimal entities in order to have an exact meaning of the scheme described above together with the principle of classical transfinite induction. However, the internal set definition principle and internal induction are a convenient tool. Clearly, both principles have their counterpart in ZFC, or, more naturally in NBG, where they can be rephrased with the transfinite induction principle. Therefore, the following argument can be rephrased in a functional analytic setting with explicit infinitesimals such as Connes' theory. In order to define the theoretical minimum for the argument, recall the transfinite induction principle. First, a nonempty linear ordered set $S$ is called well-ordered if any nonempty subset of $S$ as a least element. Next recall
\begin{defi}
An ordinal number is a transitive set which is well-ordered by the relation $\in$. The class of all ordinals is denoted by $On$. Furthermore, for $\alpha\in On$ we define $\alpha+1:=\alpha\cup \left\lbrace \alpha\right\rbrace$ to be the successor ordinal.
\end{defi}
The well-known transfinite induction principle (proofs can be found in standard text books of set theory) then can be stated as follows.
\begin{thm}
Let $R\subset On$ be a class of ordinals where $On$ is the class of all ordinals, and let $\phi$ be a property. Assume that
\begin{itemize}
 \item[i)] $\oslash\in R$ and $\phi(\oslash)$ is valid;
 \item[ii)] if $\alpha\in R$, then $\alpha +1\in R$, and if $\phi(\alpha)$ holds, then $\phi(\alpha+1)$ holds;
 \item[iii)] if $\alpha\neq \oslash$ is a nonzero limit ordinal, where $\beta \in R$ for  all $\beta \in \alpha$ and such that $\phi(\beta)$ holds for all  $\beta \in \alpha$, then $\alpha \in R$ and $\phi(\alpha)$ holds.
\end{itemize}
Then $R$ is the class of all ordinals and $\phi(\alpha)$ holds for all $\alpha \in On$. Moreover, the transfinite induction principle can be extended to every transitive class $T$ where we replace items $[i)]-[iii)]$ by just two items $a)$ and $b)$.
\begin{itemize}
 \item[a)] $\oslash\in T$ and $\phi(\oslash)$ is valid;
 \item[b)] if $\alpha\in T$ and $\phi(\beta)$ holds for all $\beta \in \alpha$, then $\phi(\alpha)$ holds.
\end{itemize}
Then for every $\alpha \in T$ $\phi(\alpha)$ holds.
\end{thm}
The latter theorem may be applied to the transitive class of ordinal numbers $T\subset On$ which is itself linearly ordered by $\in$). Transfinite induction can be used directly only for transitive classes of course, and it depends on the construction and the framework for a calculus with explicit infinitesimals whether transfinite induction can be applied to the scheme directly. If we a framework of nonstandard analysis based on Loos' Theorem as outlined in the preceding section then we have  to rely on the internal induction principle stated in Theorem \ref{intind}. Note that for our purposes even an overflow principle is sufficient, as this leads to time local regular upper bound preservation for finite positive real number distance. Finite induction and the semi-group property then lead to the time global regular upper bound. The theorem to be used for this line of argument
\begin{thm}
Let ${\mathbb U}\stackrel{*}{\rightarrow}{\mathbb U}'$ be an universe embedding as described in the preceding section, and let $\phi(x)$ be an internal $L_{{\mathbb U}'}$-formula. If $\phi(t)$ holds for all $t>t_0$ where $t-t_0$ is infinitesimal, then there exists a $t_1$ such that the shadow of $t_1-t_0$ is a positive real number $r\in {\mathbb R}_+$ such that $\phi(t_1)$ is true.   
\end{thm}
Note that related overflow principles can be used to confirm the truth of $\phi(t)$ for hyperreal values between $t_0$ and $t_1$.

However we prefer to work more closely in the framework proposed in the preceding section and consider the hyperfinite induction principle in the context of a hyperfinite set, where we consider a hyperfinite function $v^{*}_i,~1\leq i\leq n$ on a hyperfinite time interval of length $T >0$ of the form
\begin{equation}
I_{N}:=\left\lbrace  0,\frac{1}{N}T,\Delta\frac{2}{N},\cdots , \frac{N-1}{N}T,T\right\rbrace,
\end{equation} 
where $N$ is a hyperfinite number (note that hyperfinite numbers are also denoted by ${\mathbb N}^*$. 
\begin{rem}
Note that in auto-controlled schemes the step size $\delta t$ in $t$-coordinates corresponds to a step size $\frac{\Delta}{\sqrt{1-\Delta^2}}$ in transformed $s$ coordinates. As we use an equidistant time discretization here, we may write
\begin{equation}
\delta s=\frac{\Delta}{\sqrt{1-\Delta^2}}\frac{1}{N},~\delta t=\frac{\Delta}{N}.
\end{equation}
\end{rem}

Given some initial time $t_0\geq 0$   as a property to be preserved for transfinite induction we consider for $l\in I_N$ and for the finite constant $C^*>0$ of Theorem (\ref{linearboundthm}) we consider for each time step number $l$ the statement
\begin{equation}
\phi(l)\equiv {\big |}\mathbf{v}^{F}(l\delta s){\big |}^n_{h^m}\leq C^*,
\end{equation}
and the corresponding equivalent statement for the scaled function
\begin{equation}
\phi^*(l)\equiv {\big |}\mathbf{v}^{*,F}(l\delta s){\big |}^n_{h^m}\leq C^*.
\end{equation}
Note that we use the same notation $\delta t$ for the infinitesimal time step size where we tacitly assume that this is redefined accordingly if $\rho$ is different from zero. 
Here for $\mathbf{v}^{*,F}(l\delta s)=\left( v^{*,F}_1(l\delta s),\cdots ,v^{*,F}_n(l\delta s)\right)$ we use the notation 
\begin{equation}
{\big |}\mathbf{v}^{*,F}(l\delta s){\big |}^n_{h^m}:=\max_{1\leq i\leq n}{\big |}v^{*,F}_i(l\delta s){\big |}_{h^m}.
\end{equation}

The Trotter product formulas are defined with classical external sets but they can trivially extended such that the formulas in the extended universe contain only internal constants (such as ${\mathbb Z}^*$ instead of ${\mathbb Z}$ etc.). Hence, it is clear that the formulas $\phi(l)$ and $\phi^*(l)$ may be assumed to be defined by an internal formula with internal constants.
In order to apply the internal induction principle we may define the extension $\psi(l)$ for all hyperfinite numbers, where
\begin{equation}
\psi^*(l)\equiv \left\lbrace  \begin{array}{ll}
               \phi^*(l)~\mbox{ if }~l\in I_N\\
	       \\
	       \mbox{true}~\mbox{ if }~l\in {\mathbb N}^*\setminus I_N.
              \end{array}\right.
\end{equation}
We show
\begin{equation}
\forall l\in I_N :~\phi^*(l),
\end{equation}
verifies the existence of a regular time-global scheme, a statement, which corresponds to
\begin{equation}
\forall l\in {\mathbb N}^* :~\psi^*(l).
\end{equation}
The prove works for $\rho=0$. However, we mention that a similar proof holds for some $\rho\neq 0$ if certain requirements are satisfied.
If we consider an arbitrary finite time horizon $T>0$ for the original Cauchy problem for the velocity components $v_i,~1\leq i\leq n$ and the corresponding family $u^{t_0}_i,~1\leq i\leq n$ with $t_0\in [0,T)$, then this corresponds to a Cauchy problem for the scaled function $v^*_i,~1\leq i\leq n$ and the  
corresponding local time transformation family $u^{lc,*,\tau_0}_i,~1\leq i\leq n$ with $\tau_0\in [0,T_{\rho})$ with $T_{\rho}=r^{\rho}T$. Hence, in case $\rho\neq 0$ it is sufficient to sataisfy a condition of the form $r^{\rho}\geq \frac{1}{T^{\delta}}$ for some $\delta \in (0,1)$ in order to obatin a global scheme with respect to original time. Here, recall that $v^*(\tau,y)=r^{\lambda}v^{*}(r^{\rho},r^{\mu}x)=v(t,x)$. This scaling leads to the coefficients $r^{2\mu-\rho}$ for the viscosity term and the coefficient $r^{\mu-\rho+\lambda}$ for the Burgers- and Leray projection term. 

Now we choose $\rho=0$ for simplicity.
In the following we always tacitly assume that the formulas, especially Trotter product formulas are internal formulas (which can always be achieved by replacement of the external set ${\mathbb Z}$ by the internal set ${\mathbb Z}^*$).Using the internal induction principle ( or transfinite induction in a related classical argument)  it is sufficient to verify
\begin{equation}\label{ind1}
\forall \alpha\neq 0\forall l\in I_N~\max_{1\leq i\leq n}|v^{*}_{i\alpha}(l\delta t)|\leq  \frac{C^*}{1+|\alpha|^{n+s}} \mbox{ for some $s>1$}.
\end{equation}
Here for the first quantifier in (\ref{ind1}) recall our remarks on controlled schemes above which enforce the zero modes to be zero.  We consider an arbitrary fixed $s>1$ in the following. 
Note that we have for all $1\leq i\leq n$ and all modes $\alpha$
\begin{equation}\label{inda}
|v^{*}_{i\alpha}(0)|=|h_{i\alpha}(t_0)|\leq  \frac{C^*}{1+|\alpha|^{n+s}}.
\end{equation}
by assumption. Next assume inductively that for all $1\leq i\leq n$ and all modes $\alpha$ and all time step number $l\in N^*$ we have
\begin{equation}\label{indass}
|v^{*}_{i\alpha}(l\delta t)|\leq  \frac{C^*}{1+|\alpha|^{n+s}}.
\end{equation}

The at time step $l+1$ and using parameters with $\mu>0$ and $\rho=\lambda=0$ the dynamics for the modes $v^{*}_{i\alpha},~1\leq i\leq n,~\alpha\neq 0$ is given by
\begin{equation}\label{trotterlambdalaa**}
\begin{array}{ll}
v^{*}_{i\alpha}((l+1)\delta t)\doteq  {\Big (}v^{*}_{i\alpha}(l\delta t)\left(1-r^{2\mu}\nu 4\pi^2 \sum_{j=1}^n \alpha_j^2\delta t\right)\\
\\
 -r^{\mu}2\pi i \sum_{j=1}^n\sum_{\gamma \in {\mathbb Z}^n}(\alpha_j-\gamma_j)v^{*}_{i(\alpha-\gamma)}(l\delta t)v^{*}_{j\gamma}(l\delta t)\delta t\\
\\
+ r^{\mu}\frac{2\pi i\alpha_i1_{\left\lbrace \alpha\neq 0\right\rbrace}
\sum_{j=1}^n\sum_{\gamma \in {\mathbb Z}^n}\sum_{m=1}^n4\pi^2\gamma_j(\alpha_m-\gamma_m)v^{*}_{m(\alpha-\gamma)}(l\delta t)v^{*}_{j\gamma}(l\delta t)}{\sum_{i=1}^n4\pi^2\alpha_i^2}\delta t{\Big )}.
\end{array}
\end{equation}

Using the induction assumption in (\ref{indass}) and elliptic integral upper bounds of the nonlinear terms considered above we get for all $1\leq i\leq n$ and $\alpha\neq 0$
\begin{equation}\label{trotterlambdalaa**b}
\begin{array}{ll}
{\big |}v^{*}_{i\alpha}((l+1)\delta t){\big |}\lessdot  {\Big |}{\Big (}v^{*}_{i\alpha}(l\delta t)\left(1-r^{2\mu}\nu 4\pi^2 \sum_{j=1}^n \alpha_j^2\delta t\right)\\
\\
 +r^{\mu}\frac{2\pi(n+n^2)c(C^*)^2}{1+|\alpha|^{n+s}}\delta t  {\Big |}.
\end{array}
\end{equation}
for some finite constant $c$ which depends only on dimension $n$.
Next for $r>1$ we choose $\mu\geq 1$ large enough such that 
\begin{equation}\label{rmu1}
r^{\mu}=\frac{2\pi(n+n^2)c(C^*)^2}{\min\{\nu,1\}},~\mbox{ where $C^*>1$ w.lo.g.}
\end{equation}
Then for anny $1\leq i\leq n$ and any $\alpha\neq 0$ we have either 
\begin{equation}
\mbox{a)}~|v^{*}_{i\alpha}(l\delta t)|\leq \frac{C*}{2}\mbox{ or b)}|v^{*}_{i\alpha}(l\delta t)|\in \left[ \frac{C*}{2},C^*\right].
\end{equation}
In case of a) we clearly have
\begin{equation}
|v^{*}_{i\alpha}((l+1)\delta t)|\leq C^*.
\end{equation}
In case of b) we use (\ref{trotterlambdalaa**b}) and our choice in (\ref{rmu1}) in order to conclude that for all $1\leq i\leq n$ and $\alpha\neq 0$ in case of b) we have
\begin{equation}\label{trotterlambdalaa**b2}
\begin{array}{ll}
{\big |}v^{*}_{i\alpha}((l+1)\delta t){\big |}\lessdot  
{\Big |}{\Big (}C^*-\frac{C^*}{2}r^{2\mu}\nu 4\pi^2 \sum_{j=1}^n \alpha_j^2\delta t
+r^{\mu}\frac{2\pi(n+n^2)c(C^*)^2}{1+|\alpha|^{n+s}}\delta t  {\Big |}\\
  \\
{\Big |}{\Big (}C^*-\frac{C^*}{2}\frac{(2\pi(n+n^2)c(C^*)^2)^2}{\min\{\nu,1\}} 4\pi^2 \sum_{j=1}^n \alpha_j^2\delta t
 +\frac{(2\pi(n+n^2)c(C^*)^2)^2}{\min\{\nu,1\}(1+|\alpha|^{n+s})}\delta t  {\Big |}\\
 \\
 \leq \frac{C^*}{1+|\alpha|^{n+s}}
\end{array}
\end{equation}
Hence for all $1\leq i\leq n$ and $\alpha\neq 0$ we have
\begin{equation}
|v^{*}_{i\alpha}((l+1)\delta t)|\leq \frac{C^*}{1+|\alpha|^{n+s}},
\end{equation}
closing the nduction.
We note that the latter formulas are internal formulas (the quantifier ith respect to $\alpha\neq 0$ is in ${\mathbb Z}^*$) such that the internal induction principle can be implemented. Furthermore note that
we have enforced the zero modes to zero. The upper bound for the uncontrolled velocity components then gets an additional nonlinear Euler increment based on the Burgers term at each time step. We have
\begin{equation}\label{0modes}
\begin{array}{ll}
\sum_{\gamma\neq 0}e^{*}_{ij\alpha\gamma}(l \delta t)v^{*}_{i\gamma}((l\delta t)|_{\alpha=0}\delta t\\
\\
=-\sum_{\gamma\neq 0}r^{\mu}2\pi i 
(-\gamma_j)v^{*}_{i(-\gamma)}(l\delta t)v^{*}_{i\gamma}(l\delta t)\delta t.
\end{array}
\end{equation}
Using the choice in (\ref{rmu1}) and the elliptic integral upper bound again we get an upper bound for this additional increment 
\begin{equation}\label{0modes}
\begin{array}{ll}
{\big |}-\sum_{\gamma\neq 0}r^{\mu}2\pi i 
(-\gamma_j)v^{*}_{i(-\gamma)}(l\delta t)v^{*}_{i\gamma}(l\delta t)\delta t{\big |}\\
\\
\leq \frac{2\pi(n+n^2)c(C^*)^2}{\min\{\nu,1\}}2\pi(n+n^2)c(C^*)^2\leq r^{2\mu}
\end{array}
\end{equation}
with the choice in (\ref{rmu1}). Hence the linear time regular upper bound in the statement of theorem \ref{linearboundthm} follows.

\end{document}